\documentclass[twoside,leqno]{article}

\usepackage[letterpaper]{geometry}

\usepackage{siamproceedings}

\usepackage[T1]{fontenc}
\usepackage{amsfonts}
\usepackage{graphicx}
\usepackage{epstopdf}
\usepackage{enumitem}
\usepackage{color}
\usepackage{algorithmic}

\usepackage{amsmath,amssymb,graphics}
\usepackage{hyperref}
 \usepackage{comment}
 \usepackage{tabularx}
 \usepackage[protrusion=true,expansion=true]{microtype}
\usepackage{enumerate}
 \usepackage{bbm}
\usepackage{mathrsfs}
 \usepackage{paralist}

\ifpdf
  \DeclareGraphicsExtensions{.eps,.pdf,.png,.jpg}
\else
  \DeclareGraphicsExtensions{.eps}
\fi

\newsiamremark{remark}{Remark}
\newsiamremark{conjecture}{Conjecture}
\newsiamremark{hypothesis}{Hypothesis}
\crefname{hypothesis}{Hypothesis}{Hypotheses}
\newsiamthm{claim}{Claim}
\newcommand{\C}{\mathbbm{C}}

\newcommand{\N}{\mathbbm{N}}

\newcommand{\Z}{\mathbbm{Z}}
\newcommand{\R}{\mathbbm{R}}
\renewcommand{\P}{\mathbbm{P}}
\newcommand{\bbH}{\mathbbm{H}}

\newcommand{\eps}{\varepsilon}

\DeclareMathOperator{\SLE}{SLE}

\def\cV{\mathcal{V}}

\def\cS{\mathcal{S}}
\def\cR{\mathcal{R}}

\def\cL{\mathcal{L}}

\def\cE{\mathcal{E}}

\def\cA{\mathcal{A}}

\newcommand{\aryb}{\begin{eqnarray*}}
\newcommand{\arye}{\end{eqnarray*}}
\def\alb#1\ale{\begin{align*}#1\end{align*}}
\newcommand{\eqb}{\begin{equation}}
\newcommand{\eqe}{\end{equation}}
\newcommand{\eqbn}{\begin{equation*}}
\newcommand{\eqen}{\end{equation*}}

\newcommand{\BB}{\mathbbm}
\newcommand{\ol}{\overline}

\newcommand{\op}{\operatorname}

\newcommand{\frk}{\mathfrak}

\newcommand{\rta}{\rightarrow}

\newcommand{\wt}{\widetilde}
\newcommand{\wh}{\widehat} 
\newcommand{\mcl}{\mathcal}

\newcommand{\tree}{\frk T}

\newcommand{\cc}{ {\bf c}}
\newcommand{\ccL}{ {\bf c}_{\op{L}} }
\newcommand{\ccM}{ {\bf c}_{\op{M}} }

\usepackage{amsopn}

\begin{document}

\title{\Large Liouville quantum gravity: from random planar maps to conformal field theory}
    \author{Nina Holden\thanks{Courant Institute of Mathematical Sciences, New York University (\email{nina.holden@nyu.edu}).}
    \and Xin Sun\thanks{Beijing International Center for Mathematical Research, Peking University, Beijing, China
  (\email{xinsun@bicmr.pku.edu.cn}).}}

\date{}

\maketitle

% Copyright Statement
% When submitting your final paper to a SIAM proceedings, it is requested that you include
% the appropriate copyright in the footer of the paper.  The copyright added should be
% consistent with the copyright selected on the copyright form submitted with the paper.
% Please note that "20XX" should be changed to the year of the meeting.

% Default Copyright Statement
\fancyfoot[R]{\scriptsize{Copyright \textcopyright\ 2026 by SIAM\\
Unauthorized reproduction of this article is prohibited}}

% Depending on which copyright you agree to when you sign the copyright form, the copyright
% can be changed to one of the following after commenting out the default copyright statement
% above.

%\fancyfoot[R]{\scriptsize{Copyright \textcopyright\ 20XX\\
%Copyright for this paper is retained by authors}}

%\fancyfoot[R]{\scriptsize{Copyright \textcopyright\ 20XX\\
%Copyright retained by principal author's organization}}

%\pagenumbering{arabic}
%\setcounter{page}{1}%Leave this line commented out.

\begin{abstract}
Originating in theoretical physics, Liouville quantum gravity (LQG) has been an important topic in probability  theory and mathematical physics in the past two decades. In this proceeding, we review two aspects of this topic. The first is that LQG describes the random conformal geometry of the scaling limit of random planar maps. We highlight the convergence of random planar maps under discrete conformal embedding, where couplings between LQG and the Schramm-Loewner evolution (SLE) play a key role. The second aspect is the connection to conformal field theory (CFT). Here we highlight the interplay between Liouville CFT and the SLE/LQG coupling, the CFT description of 2D quantum gravity coupled with conformal matter, and applications to SLE and 2D statistical physics. We conclude with several open questions and future directions.
\end{abstract}

\section{Introduction.} Liouville quantum gravity (LQG) is a theory of random surfaces which originated in theoretical physics during the 1980s. In this century, it has evolved into a focal point in probability and mathematical physics, connecting a wide range of subjects including Gaussian free field, Gaussian multiplicative chaos, Schramm-Loewner evolution, random planar maps, and conformal field theory. 
We start by recalling the theory's physical motivation and summarizing the mathematical developments, before describing our specific focus.

%Finally, Section 1.3 details the specific focus of this proceeding, which concerns the connection to random planar maps and conformal field theory.

Two-dimensional (2D) quantum gravity is about understanding 
random surfaces. The physical origin of LQG is Polyakov's seminal work~\cite{polyakov-qg1}, which provides a field theoretic description of a class of random surfaces corresponding to quantum gravity coupled with conformal matters. In this framework, random surfaces are conformally embedded onto fixed 2D Riemannian manifolds, and the random geometry is encoded by the conformal factor. In particular, the area measure is of the form $e^{\gamma \phi} d^2z$, where $\gamma\in (0,2)$ is a parameter and $\phi$ is a random field that locally behaves as the Gaussian free field (GFF) while the global law is governed by a conformal field theory (CFT)  now known as the Liouville CFT. 
The conformal matter can also be described by a CFT. Each CFT has a key parameter called the \emph{central charge}, and Liouville CFT has central charge $\ccL$ which is related to $\gamma$ and the matter central charge $\ccM$ by 
\eqb
\ccL=1+6\Big(\frac{\gamma}{2}+\frac{2}{\gamma}\Big)^2,\qquad \ccM+\ccL=26.
\label{eq:c-relations}
\eqe
In particular, in the absence of conformal matter ($\ccM=0$) one has $\gamma=\sqrt{8/3}$. This is the case of pure 2D quantum gravity,  where the random surface is ``uniformly sampled''. 
See Section~\ref{subsec:4.3} for further explanations of Polyakov's picture.

Another natural approach to random surfaces and quantum gravity is through discretization. We use random planar maps such as random triangulations as a discrete model for 2D quantum gravity. Moreover, statistical physics models at criticality can be used as discrete models for conformal matters. The pure 2D quantum gravity can simply be modeled by uniform triangulations of a certain topological surface.  The presence of conformal matter would lead to non-uniform planar map models, where the weight of an individual map is proportional to the partition function of the corresponding statistical physics model.
 It is natural to conjecture that these two apparently different approaches to 2D quantum gravity are  equivalent. Namely, LQG describes the conformal geometry of the scaling limit of random planar maps weighted by   statistical physics models. 
 
A strong support for this conjecture in physics comes from the Knitznik-Polyakov-Zamolodchikov (KPZ) relation~\cite{kpz-scaling}. 
The is a quadratic relation between the scaling dimensions of conformal matters in the Euclidean geometry and their quantum counterpart,
which is derived from Polyakov's CFT picture. There are a number of cases where explicit calculations 
 were carried out in both the field theoretic and the discrete  approaches, and the results match. 
When the quantum dimensions are simpler to compute, the KPZ relation can give the Euclidean scaling dimensions for 2D statistical physics models as an application. See~\cite{shef-kpz,dup-icm} for references on the KPZ relation and Section~\ref{subsec:4.4} for further discussions.
 
Over the last two decades the probability community has built a solid and rich foundation for LQG. Key developments include: the construction of the random area and boundary length  measures for LQG surfaces based on Gaussian multiplicative chaos;  the construction of the LQG metric; the Brownian map as the metric-measure scaling limit of uniform random planar maps; the mating-of-trees framework for LQG coupled with Schramm-Loewner evolution (SLE) and their relation to scaling limits of random planar maps coupled with statistical physics models; the equivalence between the Brownian map and the LQG surfaces with $\gamma=\sqrt{8/3}$; the convergence to LQG surfaces of random planar maps under conformal embeddings;  the construction and resolution of Liouville CFT; the interaction between  the mating-of-trees framework and Liouville CFT; the rigorous formulation of Polykov's theory for quantum gravity on non-simply connected surfaces; the application of quantum gravity ideas such as the KPZ relation to 2D statistical physics. We do not attempt to survey the vast literature on these topics, but refer instead to previous and concurrent ICM proceedings \cite{schramm-icm,shef-icm,miller-icm,dup-icm,ddg-icm,legall-icm,rv-icm,werner-icm,smirnov-icm06,smirnov-icm10} and other texts, including e.g.\ \cite{ghs-mating-survey,berestycki-powell-book,lawler-book,curien-peeling-notes,legall-miermont-notes,gkr-review}. 

As compared to the other surveys, we focus on the following two topics. 
First, in Section \ref{sec:rpm} we present theorems and conjectures which say that random planar maps converge to LQG in the scaling limit, focusing in particular on convergence under conformal embedding and convergence in the mating-of-trees framework.   
Second, in Section~\ref{sec:cft} we explain the connection to CFT, including the relationship between Liouville CFT and the mating-of-trees framework, the matter-Liouville-ghost CFT description of 2D quantum gravity, and applications to 2D statistical physics models. We also give a precise mathematical definition of LQG surfaces in Section \ref{sec:lqg-surface} and present a selection of future directions in Section \ref{sec:outlook}.

\section{Liouville quantum gravity: a  theory of random surfaces.}
\label{sec:lqg-surface}

In this section we give a precise mathematical definition of LQG as a model for a random surface, following \cite{shef-kpz} and subsequent papers. 

Let $D\subset\C$ be a domain. The 2D (Dirichlet) Gaussian free field (GFF) in $D$ is the centered Gaussian field $h$ with covariance given by the zero-boundary Green's function $G_D$ of $-\frac{1}{2\pi}\Delta$, where $\Delta$ is the Laplace operator on $D$. For fixed $x\in D$ we have $G_D(x,y)=\log |x-y|+c(x)+o(1)$ as $y\rta x$ for some constant $c(x)$ depending only on $x$. In particular, $G_D(x,x)=\infty$, which means that $h$ cannot be well-defined as a function. One can argue that instead $h$ is well-defined as a random distribution, i.e., a random generalized function. In fact, $h$ is a random variable almost surely taking its values in the Sobolev space $H^{-\eps}(D)$ for any $\eps>0$. The GFF arises as the scaling limit of a number of functions appearing in statistical physics and it plays an important role in modern probability theory and theoretical physics.

Let $\gamma\in(0,2)$. Heuristically speaking, a $\gamma$-LQG surface is a 2D Riemannian manifold with metric $e^{\gamma h}(dx^2+dy^2)$, where $h$ is a GFF or some closely related random distribution and $dx^2+dy^2$ is the standard Euclidean metric. This definition of an LQG surface does not make rigorous sense since the GFF is a distribution and not a function. However, several observables of the surface, such as its area measure and metric, can be defined in a rigorous way via regularization of the field. 

The area measure can be defined rigorously via the theory of Gaussian multiplicative chaos (GMC), which was studied in \cite{kahane,hk-gmc} and more recently in e.g.\ \cite{shef-kpz,rhodes-vargas-review,robert-vargas-revisited,berestycki-gmt-elementary,berestycki-powell-book,shamov-gmc,shef-wang-lqg-coord}. For any $z\in D$ and $\epsilon>0$ such that the ball $B_\epsilon(z)$ is contained in $D$, the average $h_\epsilon(z)$ of $h$ on the circle $\partial B_\eps(z)$ is well-defined. Furthermore, for any $\eps>0$ the map $z\mapsto h_\epsilon(z)$ has a continuous modification. Therefore the measure $\cA_\epsilon$ given by $\cA_\epsilon(d^2z)=\eps^{\gamma^2/2}e^{\gamma h_\epsilon(z)}\,d^2z$ is well-defined. It can be argued that this measure converges a.s.\ as $\eps\rta 0$ to some limiting measure which we denote by $\cA$ or $\cA_h$. It can also be argued that a number of other approximation schemes give the same limiting measure $\cA_h$. Furthermore, one can define an LQG boundary measure $\cL_h$ supported on $\partial D$ in a similar manner as the LQG area measure $\cA_h$.

The LQG metric can also be constructed via regularization of the field, but the construction is much harder than for the area measure and it was only completed recently in \cite{dddf-lfpp,gm-uniqueness}. Ding, Dubedat, Dunlap, and Falconet \cite{dddf-lfpp} studied metrics associated to a GFF mollified by the two-dimensional heat kernel. They proved that suitably renormalized versions of this metric is tight with respect to the uniform topology. In \cite{dfgps-metric} it was proven that any subsequential limiting metric satisfies a certain list of natural axioms. Finally, Gwynne and Miller \cite{gm-uniqueness} proved uniqueness of subsequential limits and that the unique subsequential limit is a strong LQG metric $\frk D_h:D\times D\to [0,\infty)$, which is a metric satisfying locality, Weyl scaling, and a coordinate change formula for translation and scaling, such that the resulting metric space is a length space. The Hausdorff dimension $d_\gamma$ of a $\gamma$-LQG surface is not known for general $\gamma$, although close upper and lower bounds have been established \cite{dg-lqg-dim,ang-discrete-lfpp,gp-lfpp-bounds} and it is known that the function $\gamma\mapsto d_\gamma$ is strictly increasing. For $\gamma=\sqrt{8/3}$ the metric can be equivalently constructed via the growth process known as the Quantum Loewner Evolution (QLE) and the resulting metric space is equal in law to a Brownian surface \cite{lqg-tbm3}.

We view an LQG surface as an abstract surface which can be embedded conformally into the complex plane in several different ways. Each embedding can be associated with pair $(D,h)$, where $D\subset \C$ is a domain and $h$ is the field used to define the measure and the metric. We say that two pairs $(D,h)$ and $(\widetilde D,\widetilde h)$ are equivalent if there is a conformal transformation $\phi:D\to\wt D$ such that $\cA_h(A)=\cA_{\wt h}(\phi(A))$ for every measurable $A\subset D$ and $\frk D_h(x,y)=\frk D_{\wt h}(\phi(x),\phi(y))$ for every $x,y\in D$. This condition is equivalent to $h=\wt h\circ\phi+Q\log|\phi'|$, where 
$$Q=2/\gamma+\gamma/2$$
is the background charge. The $\gamma$-LQG surface (i.e., the equivalence class) associated with the pair $(D,h)$ is denoted by $(D,h)/{\sim_\gamma}$. One can also consider LQG surfaces $(D,h,z_1,\dots,z_n,\eta_1,\dots,\eta_m)$ with marked points $z_1,\dots,z_n\in\ol D$ and/or curves $\eta_1,\dots,\eta_m$ for $n,m\in\N$, where we say that two such tuples are equivalent if the conformal map $\phi$ transforms the points $z_1,\dots,z_n$ and the curves $\eta_1,\dots,\eta_m$ covariantly.

The measure and metric discussed above can be defined for several fields $h$ closely related to the GFF, for example for a field which is absolutely continuous with respect to a GFF plus a (possibly random) bounded continuous function. 
As we will see in Sections \ref{sec:rpm} and \ref{sec:cft}, among all such variants of the GFF there are certain choices for $h$ that are particularly natural since they arise as the scaling limit of random planar maps and/or appear in LCFT. %As we will see in Section \ref{} below, in several contexts it has been proven that a field which describes the scaling limit of random planar maps is identical to the corresponding field in LCFT. 

One example is the surface known as the \emph{quantum disk}, with law denoted by $\op{QD}$, which is the canonical surface with disk topology. There is also a generalized version $\op{QD}_{n,m}$ of this disk with $n\in\N$ singularities in the bulk and $m$ singularities on the boundary. Conditioned on a surface sampled from $\op{QD}_{n,m}$ with the location of the singularities forgotten, the singularities have the law of points sampled uniformly from the area measure $\cA_h$ and the boundary measure $\cL_h$, respectively. For $n=0$ and $m=2$ one can also, via the mating-of-trees framework, define a variant $\op{QD}_2^W$ of this disk with generalized boundary singularities of weight $W>0$. Nearby a singularity $z_0$ on the boundary of weight $W$ the field is absolutely continuous with respect to a free boundary GFF plus $\beta\log|\cdot-z_0|^{-1}$, where 
$\beta = \gamma+2/\gamma-W/\gamma$. 
The case $W=2$ plays a special role since the field looks like a typical point sampled from the LQG boundary measure $\cL_h$ nearby the singularities; more precisely, we have $\op{QD}_{0,2}=\op{QD}^W_2$ for $W=2$.  
Similarly, the canonical measure on LQG surfaces with spherical topology is called the \emph{quantum sphere} and the different variants are denoted by $\op{QS}$, $\op{QS}_n$, and $\op{QS}_2^W$ for $n\in \N$ and $W>0$. 

%There is also a generalized version $\op{QD}_n^W$ of this disk with $n\in\N$ singularities on the boundary of weight $W>0$. Nearby a singularity $z_0$ of weight $W$ the field is absolutely continuous with respect to a free boundary GFF plus $\beta\log|\cdot-z_0|$, where 
%$$
%\beta = \gamma+2/\gamma-W/\gamma.
%$$ 
%The case $W=2$ plays a special role since the field looks like a typical point sampled from the LQG boundary measure $\cL_h$ nearby the singularities and we therefore skip $W$ from the notation in that case, i.e., we have $\op{QD}_n=\op{QD}_n^W$ with $W=2$. Similarly, the canonical measures on LQG surfaces with spherical topology is called the \emph{quantum sphere} and the different variants are denoted by $\op{QS}$, $\op{QS}_n$, and $\op{QS}_n^W$, where now $\op{QS}_n=\op{QS}_n^W$ for $W=4-\gamma^2$.

There are canonical ways to disintegrate the LQG disk  measures according to the LQG length of the various boundary segments. For example, there are measures $\op{QD}_2^W(\ell_{\op{L}},\ell_{\op{R}})$ for 
$W,\ell_{\op{L}},\ell_{\op{R}}>0$ 
such that the boundary segments between the two marked points have LQG lengths $\ell_{\op{L}}$ and $\ell_{\op{R}}$, respectively, and such that 
$$
\op{QD}_2^W= \int_0^\infty\int_0^\infty
\op{QD}_2^W(\ell_{\op{L}}, \ell_{\op{R}})\,
d\ell_{\op{L}}\,d\ell_{\op{R}}.
$$
Similarly, we can define measures $\op{QD}(\ell)$ and $\op{QD}_{0,1}(\ell)$ supported on disks of boundary length $\ell>0$.
These examples will be used in Section~\ref{sec:sle-welding-mot}.

\section{LQG and random planar maps.}
\label{sec:rpm}
\begin{figure}
    \centering
    \includegraphics[scale=1]{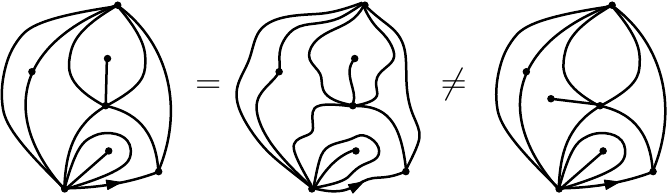}
    \caption{Illustration of the definition of a planar map. All planar maps we consider are rooted, meaning that they have a distinguished oriented edge indicated by an arrow. The leftmost and the rightmost planar maps are equivalent viewed as graphs, but viewed as planar maps they are different.}
    \label{fig:planar-map-def}
\end{figure}

A planar map is a proper embedding of a connected (multi-)graph in the sphere $\mathbb S^2$, where two planar maps are viewed as equivalent if one can get one from the other by doing an orientation-preserving homeomorphism from the sphere to itself; see Figure \ref{fig:planar-map-def}. Planar maps are studied in many different branches of both mathematics (combinatorics, random geometry, random matrix theory, ...) and physics (string theory, conformal field theory, gauge theory, statistical mechanics, ...). In particular, in probability theory and theoretical physics random planar maps are viewed as models for random surfaces. Random planar maps are partly of interest since they are natural models for random surfaces and partly since they have deep applications in theoretical physics.

A research direction in the study of random planar maps is understanding their scaling limits, motivated by their role in 2D quantum gravity and as canonical models for random surfaces. Confirming predictions from the physics literature \cite{polyakov-qg1,kpz-scaling,david-conformal-gauge,dk-qg}, it has been proven that a number of natural random planar maps, including planar maps decorated by statistical physics models, converge in the scaling limit to LQG surfaces. To describe more precisely the general scaling limit conjecture, 
for $n\in\N$ suppose $M_n$ is a planar map with $n$ edges sampled at random such that for any fixed map $\frk m$ with $n$ edges,
\eqb
\P[M_n=\frk m] = c_0(\op{det}\Delta_{\frk m})^{-\ccM/2},
\label{eq:rpm}
\eqe
where $\Delta_{\frk m}$ is the Laplacian of the planar map,  and $c_0>0$ is a normalizing constant. 
The choice of weight on the right side of~\eqref{eq:rpm} ensures that we get a model for 2D  quantum gravity coupled with conformal matter that has central charge $\ccM$.
\begin{conjecture}
For $\ccM\leq 1$, the planar maps $M_n$ converge in the scaling limit to a $\gamma$-LQG surface, where $\gamma$ and $\ccM$ are related as in \eqref{eq:c-relations}.
	\label{conj1}
\end{conjecture}
Many variants of this conjecture are also believed to hold. For example, on the right side of \eqref{eq:rpm} we can consider the partition function of some critical statistical physics model on $\frk m$. A variety of 2D  statistical physics models at criticality can be used to describe conformal matters with $\ccM\le 1$. For example, the percolation model and the Ising model correspond to a conformal matter  with $\ccM=0$ and $\ccM=\frac12$, respectively.   
% the partition function at criticality is believed to be a good approximation to the right side of \eqref{eq:rpm} for $\ccM$ equal to the central charge of the model. 
Yet another alternative is to require that $M_n$ belongs to some subset of the set of planar maps with $n$ edges, e.g.\ the set of triangulations (which are planar maps where all faces are required to have degree three) or simple maps (which disallow multiple edges and loops), or to require that the planar map has disk topology (which means, roughly speaking, that it has a macroscopic boundary, sometimes required to be a simple curve). 

The conjecture above is stated quite imprecisely since we did not explain what it means for random planar maps to converge to LQG. In fact, one can use several different topologies when defining convergence of planar maps to LQG. There are three commonly used notions: convergence as a metric measure space (Gromov-Hausdorff-Prokhorov topology), convergence in the mating-of-trees framework, and convergence under conformal embedding. We focus on the last two notions of convergence in this survey. 

In Section \ref{sec:sle-welding-mot} we give relevant background on SLE, conformal welding, and mating of trees, which will be used throughout the remainder of this survey. In Section \ref{sec:discrete-mot} we present mating-of-trees bijections and some applications, such as convergence of planar maps in the mating-of-trees framework. Finally, in Section \ref{sec:conformal-embedding} we present our result on convergence of uniform triangulations to $\sqrt{8/3}$-LQG under a discrete conformal embedding called the Cardy-Smirnov embedding.

\subsection{The Schramm-Loewner evolution, conformal welding, and mating of trees.}
\label{sec:sle-welding-mot}
The Schramm-Loewner evolution (SLE) is a one-parameter family of random fractal curves. They were introduced by Schramm \cite{schramm0} as a candidate for the scaling limit of various 2D statistical physics models at criticality and in later works it was proven that the SLE indeed describes such limits \cite{lsw-lerw-ust,smirnov-cardy,smirnov-ising}. The most classical variant of SLE is the chordal SLE$_\kappa$ for $\kappa>0$, which is a curve connecting two given boundary points of a simply connected domain. Chordal SLE$_\kappa$ is uniquely characterized by two natural properties known as conformal invariance and the domain Markov property. A more explicit definition of the chordal SLE$_\kappa$ can be obtained by solving the Loewner differential equation with driving function equal to $\sqrt{\kappa}$ times a standard Brownian motion. 

Properties of the SLE$_\kappa$ depend strongly on $\kappa$ \cite{schramm-sle}. For $\kappa\in(0,4]$ the chordal SLE$_\kappa$ is a simple curve. For $\kappa\in(4,8)$ the chordal SLE$_\kappa$ is a curve which hits its past and the domain boundary. For $\kappa\geq 8$ the chordal SLE$_\kappa$ is space-filling, meaning that its trace is equal to the closure on the considered domain. There is also a space-filling variant of SLE$_\kappa$ for $\kappa\in(4,8)$, which can be constructed via the theory of imaginary geometry. Imaginary geometry is an approach to SLE where SLE$_\kappa$ curves are constructed as flow lines of a vector field $e^{ih/\chi}$ for $\chi=2/\sqrt{\kappa}-\sqrt{\kappa}/2$ and $h$ equal to a GFF \cite{dubedat-coupling,ig1,ig4}. 

\begin{figure}
    \centering
    \includegraphics{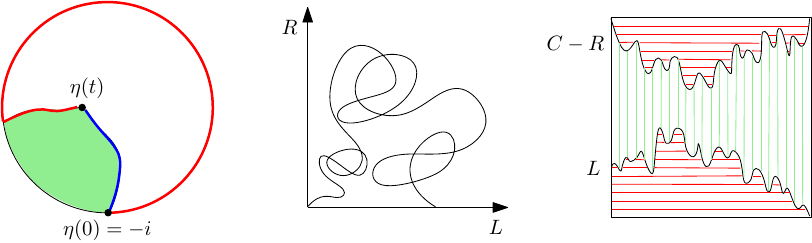}
    \caption{ 
    {\bf Left:} Illustration of a clockwise space-filling SLE$_\kappa$ loop $\eta$ in $\BB D$ starting and ending at $-i$ for $\kappa\geq 8$. The red curve has quantum length $L_t$ while the blue curve has quantum length $R_t$. 
    {\bf Middle:} The boundary length process $(L,R)$ of $\eta$ has the law of a 2D Brownian path measure from $(0,1)$ to $(0,0)$ in the first quadrant.
    {\bf Right:} The black curves are $L$ and $C-R$ for $C>0$ chosen sufficiently large such that these two curves do not intersect. By identifying points lying on the same (red) horizontal line above $C-R$ or below $L$, or on the same (green) vertical line between the two curves, we obtain the left figure. This procedure is called mating of trees.}
    \label{fig:mot}
\end{figure}

There are a number of variants of SLE$_\kappa$ beyond the chordal SLE$_\kappa$. These curves have the same local properties as the chordal SLE$_\kappa$, but the global geometry is different. A classical example is the conformal loop ensemble (CLE$_\kappa$) \cite{shef-cle,shef-werner-cle}, which is a countable collection of SLE$_\kappa$-type loops. Another example, which will be studied later in this section, is the clockwise space-filling SLE$_\kappa$ loop for $\kappa\geq 8$, which is the limit of chordal SLE$_\kappa$ as the terminal point approaches the initial point in clockwise direction; see the left panel of Figure \ref{fig:mot}. The clockwise space-filling SLE$_\kappa$ loop also exists for $\kappa\in(4,8)$ as the same limit of the chordal space-filling SLE$_\kappa$.

Besides the definition of SLE via the Loewner equation and imaginary geometry, there is a third way of generating the SLE, namely via the conformal welding of LQG surfaces. We illustrate this construction using the example of two LQG surfaces with disk topology and two boundary marked points.
%It has been proven in various settings that LQG surfaces can be glued (more precisely, conformally welded) together such that the interface of the resulting surface is given by an SLE. 
Given two LQG surfaces $(\BB D,h_{\op{L}}a_{\op{L}},b_{\op{L}})/{\sim_\gamma}$ and $(\BB D,h_{{\op{R}}},a_{\op{R}},b_{\op{R}})/{\sim_\gamma}$, each with two marked boundary points, suppose that 
the LQG length assigned to the counterclockwise boundary arc $\partial_{a_{\op{L}},b_{\op{L}}}$ from $a_{\op{L}}$ to $b_{\op{L}}$ by $h_{\op{L}}$ is equal to  
the LQG length assigned to the counterclockwise boundary arc $\partial_{b_{\op{R}},a_{\op{R}}}$ from $b_{\op{R}}$ to $a_{\op{R}}$ by $h_{\op{R}}$. 
Then we can define a homeomorphism $\phi:\partial_{[a_{\op{L}}, b_{\op{L}}]} \to \partial_{[b_{\op{R}}, a_{\op{R}}]}$ which  preserves LQG length and satisfies $\phi(a_{\op{L}})=a_{\op{R}}$ and $\phi(b_{\op{L}})=b_{\op{R}}$. A conformal welding of the two disks using the identification $\phi$ is a conformal structure on
the disk obtained by identifying 
$\partial_{[a_{\op{L}}, b_{\op{L}}]}$ 
and 
$\partial_{[b_{\op{R}}, a_{\op{R}}]}$ 
according to $\phi$. More precisely, one wants to find a (chordal) curve $\eta$ in $\BB D$ along with conformal transformations $\psi_{\op{L}}$ and $\psi_{\op{R}}$ from $\BB D$ to the two
components of $\BB D\setminus \eta$ such that $\phi=\psi_{\op{R}}^{-1}\circ\psi_{\op{L}}$. If such a curve $\eta$ and conformal maps $\psi_{\op{L}}$ and $\psi_{\op{R}}$ exist and are unique (modulo conformal maps from $\BB D$ to itself) then we say that the conformal welding of the two surfaces is well-defined.

A number of conformal welding results have been established for LQG surfaces, starting with Sheffield \cite{shef-zipper} and further developed by Duplantier, Miller, and Sheffield \cite{wedges} and later in e.g.\ \cite{ahs-disk,ahs-loop} and subsequent works. In many cases one is able to identify explicitly the law of the surface and the curve $\eta$ obtained upon welding, and also prove independence of the surface and the curve. The parameter $\gamma\in(0,2]$ of the LQG and the parameter $\kappa>0$ of the SLE are always related via $\kappa\in\{ \gamma^2,16/\gamma^2 \}$. Many of the conformal welding results are inspired by bijections for random planar maps, although the proofs are done only in the continuum.
\begin{figure}
    \centering
    \includegraphics[scale=0.95]{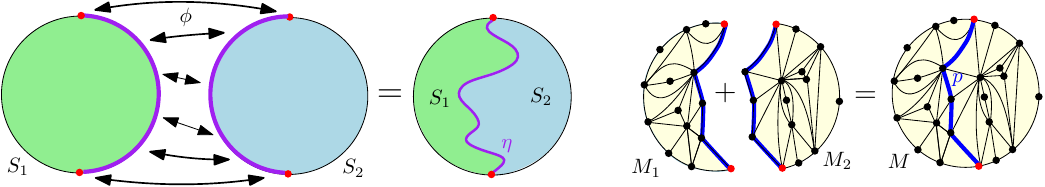}
    \caption{{\bf Left:} The conformal welding of two quantum disks $S_1$ and $S_2$ sampled independently from $\op{QD}_2$ (conditioned on the event that the two purple curves in the left part of the figure have the same quantum length) is a quantum disk sampled from $\op{QD}_2^W$ with $W=4$ decorated by an independent chordal SLE$_\kappa$ $\eta$.
    {\bf Right:} The discrete counterpart of the conformal welding result in the left figure for $\kappa=\gamma^2=8/3$. The figure is illustrating that there is a bijection between the objects on the left and right sides of the equality, where on the right side we have a planar map $M$ with disk topology and two marked boundary points, along with a self-avoiding curve $p$ connecting the two marked points.}
    \label{fig:welding}
\end{figure}

For concreteness we include the following example due to the authors and Ang \cite{ahs-disk}, where Weld$(\cdot,\cdot)$ denotes conformal welding of the quantum surfaces in question, $\op{QD}_2^W(\ell_{\op{L}},\ell_{\op{R}})$ for $W,\ell_{\op{L}},\ell_{\op{R}}>0$ is the two-pointed quantum disk with boundary lengths $\ell_{\op{L}}$ and $\ell_{\op{R}}$ and weight $W>0$ introduced in Section \ref{sec:lqg-surface}, and the notation on the left side means that an independent SLE$_\kappa$ is drawn on the quantum surface connecting its two marked points. The theorem gives an identity between two  measures on the space of curve-decorated quantum surfaces. See Figure \ref{fig:welding} for an illustration. The general result in~\cite{ahs-disk} states that the conformal welding of a sample from  $\op{QD}_2^{W_1}(\ell_{\op{L}},\ell)$  and  a sample from $\op{QD}_2^{W_2}(\ell,\ell_{\op{R}})$ gives a sample from $\op{QD}_2^{W_1+W_2}(\ell_{\op{L}},\ell_{\op{R}})$ coupled with an explicit variant of SLE. We only present the case where $W_1=W_2=2$ to avoid introducing this variant of SLE.
\begin{theorem}\label{thm1}
   There is a constant $c>0$ such that for all $\ell_{\op{L}},\ell_{\op{L}}>0$ and $\kappa=\gamma^2\in(0,4)$,
    $$\op{QD}_2^4(\ell_{\op{L}},\ell_{\op{R}}) \otimes\SLE_\kappa= c \int_0^\infty  \op{Weld}\big(\op{QD}_2^2(\ell_{\op{L}},\ell)\times \op{QD}_2^2(\ell,\ell_{\op{R}})\big)  \,d\ell.$$
\end{theorem}

The paper \cite{ahs-disk} also proves a sphere version of this theorem which says that the conformal welding of two copies of $\op{QD}_2^2$ (independent conditioned on matching boundary lengths) into an LQG surface with spherical topology is given by $\op{QS}_2^4$ decorated by an explicit variant of SLE.

Mating of trees is another powerful coupling between LQG surfaces and SLE which builds on the conformal welding results described above. To state a precise version, consider a unit boundary length quantum disk with one marked boundary point $(\BB D,h,1)/{\sim_\gamma}$ sampled from $\op{QD}_1(1)$ and draw an independent clockwise space-filling SLE$_\kappa$ loop on top of this disk starting and ending at the marked point $1$ and with $\kappa=16/\gamma^2>4$. Let $A:=\cA_h(\BB D)$ denote the total LQG area of the disk and parameterize $\eta$ by LQG area so that $\cA_h(\eta([0,t]))=t$ for any $t\in[0,A]$. The pair $(h,\eta)$ determines a stochastic  process $(L_t,R_t)_{t\in[0,A]}$ such that $L_t$ (resp.\ $R_t$) is the left (resp.\ right) boundary length of the set $\eta[0,t]$; see Figure~\ref{fig:mot}. It can be proven that $(L,R)$ has the law of a Brownian path measure in the first quadrant from $(1,0)$ to the origin. The correlation between two coordinates is given by $-\cos(4\pi/\kappa)$ \cite{wedges,kappa8-cov}. Furthermore, $(h,\eta)$ is measurable with respect to $(L,R)$ and can be obtained by ``mating'' the pair of trees with height functions $L$ and $R$.

The previous paragraph describes the disk version of mating of trees. There is also a whole-plane version, where the boundary length process is a two-sided Brownian motion and the LQG surface is a so-called quantum cone. Furthermore, there is a sphere version where the boundary length process is a Brownian path in the first quadrant starting and ending at the origin and the LQG surface is the quantum sphere $\op{QS}_1$ with one marked point. See~our survey with Gwynne on further background on mating of trees~\cite{ghs-mating-survey}.

\subsection{Discrete mating of trees and applications.}
\label{sec:discrete-mot}
There are discrete counterparts of mating of trees in the context of planar maps. In fact, such discrete results inspired the development of the continuum mating of trees. The first-discovered  example is the Mullin bijection~\cite{mullin-maps,bernardi-maps,shef-burger}, see also \cite{chen-fk,ghs-mating-survey,ghs-map-dist}. 
\begin{theorem}
    For any $n\in\N$ there is a bijection between
    \begin{compactitem}
        \item pairs $(M,\frk T)$, where $M$ is a planar map with $n$ edges and $\frk T\subset\cE(M)$ is a spanning tree on $M$, and
        \item nearest-neighbor walks $(\cL_k,\cR_k)_{0\leq k\leq 2n}$ on $\Z_+^2$ starting and ending at $(0,0)$.
    \end{compactitem}
\end{theorem}
\begin{figure}
    \centering
    \includegraphics[scale=0.8]{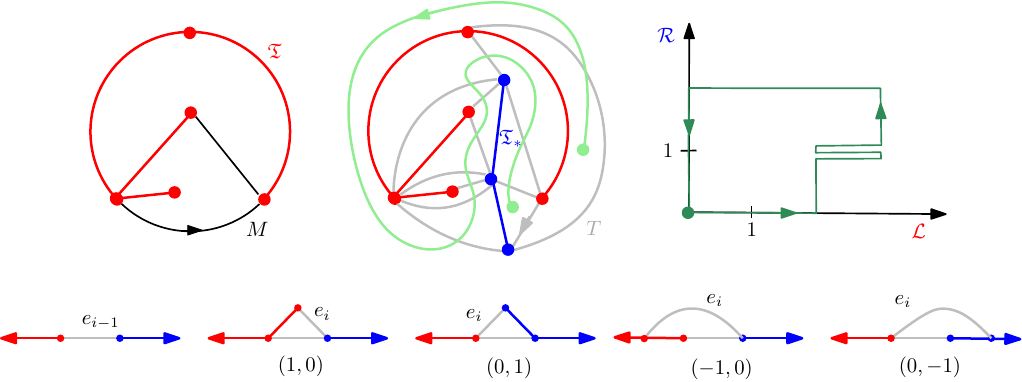}
    \caption{Illustration of the Mullin bijection. {\bf Top left}: A planar map $M$ with a spanning tree $\frk T$. 
    {\bf Top middle}: The triangulation $T$ with trees $\frk T$ and $\frk T_*$. The green curve indicates the order in which the faces (triangles) of the map were added in the sewing procedure. 
    {\bf Top right}: The walk $(\cL,\cR)$ associated with $(M,\frk T)$ via the Mullin bijection.
    {\bf Bottom}: Illustration of the sewing procedure. Before reading the $k$th step of the walk, the part of the map we constructed so far is ``below'' the line in the left panel. After reading the $k$th step of the walk the map looks like one of the other four panels.
    }
    \label{fig:mullin}
\end{figure}
Let us give a brief description of the bijection. First denote the set of vertices, edges, and faces of the planar map $M$ by $\cV(M)$ and $\cE(M)$, respectively. Let $M_*$ be the dual map of $M$ and let $\frk T_*$ be the dual spanning tree, i.e., $\frk T_*$ is the set of edges of $M_*$ which do not cross edges of $\tree$. Let $T$ be the triangulation with vertex set $\mcl V(M)\cup \mcl V(M_*)$ and edge set given by the union of $\frk T\cup \frk T_*$ and edges of the form $\{v,v_* \}$, where $v\in \cV(M)$ is incident to the face that is dual to $v_*\in \cV(M_*)$. We assign the root of $T$ in specific way such that the tuples $(M,\frk T)$ and $(T,\frk T,\frk T_*)$ determine each other. Given the walk $(\cL,\cR)$ we can construct $(T,\frk T,\frk T_*)$ by reading the steps of the walk one by one and do a procedure known as the \emph{sewing procedure}. At time $k=0$ we start from a single edge $e_0$ and each time we read another step of the walk we add a triangle as shown in Figure \ref{fig:mullin}. The $k$th step of the walk is associated with an edge $e_k$ such that the $k$the triangle we add has edges $e_{k-1}$ and $e_k$ on its boundary. When we have finished reading the steps of the walk we identify the edges $e_0$ and $e_{2n}$, which gives us $(T,\frk T,\frk T_*)$. 

When $n\rta\infty$ the process $(\cL,\cR)$ appropriately rescaled converges in law to a 2D Brownian excursion $(L,R)$ in the first quadrant starting and ending at the origin. We interpret this as a convergence result of the tree-weighted planar map to an SLE-decorated quantum sphere, where the (green) Peano curve in the middle panel of Figure \ref{fig:mullin} is the discrete counterpart of the space-filling SLE$_\kappa$. Since the two coordinates of the Brownion excursion $(L,R)$ are uncorrelated we have $\gamma=4/\sqrt{\kappa}=\sqrt{2}$. 

A number of other mating-of-trees bijections have also been discovered \cite{shef-burger,bhs-site-perc,kmsw-bipolar,gkmw-burger,lsw-schnyder-wood}. This includes in particular the following result from \cite{bhs-site-perc} (building on \cite{bernardi-dfs-bijection}) obtained in collaboration with Bernardi. Recall that a planar map is called a triangulation of a disk if all faces have exactly three edges except for the exterior face (i.e., the face to the right of the root edge), which has arbitrary degree (called the perimeter) and simple boundary. We consider type II triangulations, which are allowed to have multiple edges between two vertices, but which are not allowed to have loops.
%%% see Cor 2.9 in joint convergence paper
\begin{theorem}
    For any $\ell\in\Z_+$ there is a bijection between
    \begin{compactitem}
        \item pairs $(M,\omega)$, where $M$ is a triangulation of a disk with perimeter $\ell+2$ and $\omega$ is a site percolation $\omega:\cV(M)\to\{\op{blue},\op{yellow} \}$ with monochromatic blue boundary condition, and 
        \item walks in the first quadrant $\BB Z_+^2$ with steps in $\{ (0,1), (1,0), (-1,-1) \}$ starting at $(\ell,0)$ and ending at $(0,0)$.
    \end{compactitem}
    \label{thm:perc}
\end{theorem}
The percolated triangulation $(M,\omega)$ can be constructed from the walk $(\cL,\cR)$ by a sewing procedure in a similar way as for the Mullin bijection, but now $(0,1)$ and $(1,0)$ correspond to adding a new triangle while $(-1,-1)$ corresponds to identifying two edges of the already constructed map. The coordinates of the walk have correlation $1/2$, which means that we have $\gamma=4/\sqrt{\kappa}=\sqrt{8/3}$ in the scaling limit. We refer to \cite{bhs-site-perc} for further details.

The discrete and continuum variants of mating of trees have had numerous applications. First, there are purely continuum applications where information about LQG and SLE is extracted from the Brownian motion $(L,R)$ appearing in continuum mating of trees. Second, there are applications to random planar maps. Such applications include the convergence results for random planar maps to LQG provided by the mating-of-trees bijections, as we discussed above for tree-weighted planar maps and uniform triangulations. Furthermore, by considering strong couplings of $(\cL,\cR)$ and $(L,R)$ one is sometimes able to extract more information about the planar map from LQG and SLE. %By considering couplings between random walk and 2D Brownian motion we get a coupling of planar maps and LQG, which is beneficial for the study of both. 
We will describe one concrete application of the latter type from our paper with Gwynne \cite{ghs-map-dist} and refer to Section \ref{sec:conformal-embedding} and \cite{gwynne-miller-saw,lqg-tbm3,gp-dla,ghm-kpz,gm-spec-dim,gh-displacement} for further applications.

Recall from Section \ref{sec:lqg-surface} that the LQG metric can be approximated by considering a particular mollification of the GFF by the 2D heat kernel. There are several other natural (discrete or smooth) approximations to the LQG metric for which convergence to the LQG metric has not been proven. One such approximation is the random planar map known as the mated-CRT map, which is defined in terms of the planar Brownian motion $(L,R)$ arising in mating of trees.  
By using mating-of-trees bijecions for planar maps along with strong coupling of the associated random walks $(\cL,\cR)$ and $(L,R)$, we show  \cite{ghs-map-dist} that the ball growth exponent in several natural planar maps is equal to the ball growth exponent of the mated-CRT map. In particular, this was carried out for tree-weighted planar maps, bipolar-oriented planar maps, and Schnyder wood-decorated planar maps. It was later proven by Ding and Gwynne \cite{dg-lqg-dim} that this ball growth exponent is equal to the LQG dimension $d_\gamma$. This means that any bound one gets on the LQG dimension $d_\gamma$ \cite{dg-lqg-dim,ang-discrete-lfpp,gp-lfpp-bounds} automatically gives a bound for the ball growth exponent of the associated planar map model treated in \cite{ghs-map-dist}.

%and, although convergence of these approximations has not been proven, it is known that the ball growth exponent of several of these approximations is equal to the dimension $d_\gamma$ of a $\gamma$-LQG surface \cite{dg-lqg-dim}.

\subsection{Convergence of planar maps under conformal embedding.}
\label{sec:conformal-embedding}

In this section we consider convergence of planar maps under conformal embedding. This notion of convergence is closest to the spirit that LQG describes the random conformal geometry of surfaces in 2D quantum gravity.
One wants to embed the planar map into the complex plane via a so-called discrete conformal map, which is an embedding $\phi:\cV(M)\to\C$ that can be viewed as a discrete approximation of a conformal map. Some classical examples of discrete conformal embeddings are circle packing \cite{stephenson-book,RodinSullivan87,legall-icm}, embedding via the uniformization theorem for Riemann surfaces \cite{shef-kpz,curien-glimpse}, and the Tutte embedding \cite{gms-tutte}. When embedding the planar map one gets an area measure and a metric (distance function) in the plane, given by counting measure on the vertices and the graph metric, respectively; see \eqref{eq:meas-emb} and \eqref{eq:metric-emb} below for more concrete definitions. The goal is to show that these converge upon appropriate renormalization to the Liouville area measure and the Liouville metric, respectively. Alternatively, one can consider only the measure (not the metric) since metric properties of non-uniform planar maps are currently quite poorly understood. The following is a variant of Conjecture \ref{conj1} in this context.  

\begin{conjecture}
	For $n\in\N$ sample a random planar map $M_n$ as in \eqref{eq:rpm} or as discussed right below Conjecture~\ref{conj1}. Embedd the planar map into $\BB S^2$ (or a subset of $\BB C$) via a discrete conformal embedding. Then the area measure and metric induced on $\BB S^2$ (or $\BB C$) converge in the scaling limit to the area measure and metric, respectively, associated with $\gamma$-LQG. 
	\label{conj2}
\end{conjecture}

%\begin{wrapfigure}{r}{0.4\textwidth}
%	\centering
%	\vspace{-10pt}
%	\includegraphics[scale=1.1]{../../../A - Talks/Bernardi-Sun-Sitepercolation/figures/tutte-rw3}
%	\vspace{-10pt}
%	\caption{Planar map with disk topology embedded in $\BB D$ with the Tutte embedding. The red path is the trace of a simple random walk started from the blue vertex and run until hitting the boundary. Figure made by J.\ Miller.}
%	\vspace{-10pt}
%	\label{fig:tutte}
%\end{wrapfigure}

In \cite{hs-cardy-embedding}, building on  \cite{aasw-type2,ghss-ldp,hlls-cut-pts,hls-sle6,bhs-site-perc,ghs-metric-peano}, the conjecture is proven in the special case of uniform trianguations (i.e., $\ccM=0$) and $\gamma=\sqrt{8/3}$. More precisely, we consider the \emph{Boltzman measure} on the set of triangulations of a disk with given perimeter  $\ell\in\{2,3,\dots \}$. This measure gives mass proportional to $(2/27)^n$ to each fixed type II triangulation of a disk with perimeter $\ell$ and $n$ vertices. %Our proof allows both type I, type II, and type III triangulations.

\begin{figure}
    \centering
	\includegraphics[scale=1]{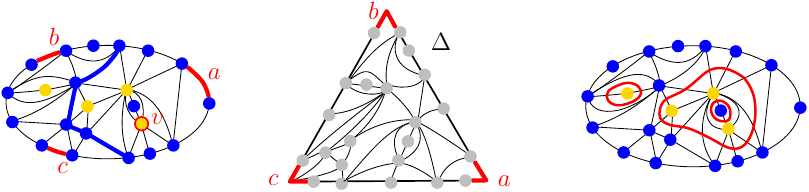}
	\caption{
    {\bf Left}: Illustration of the event the $E_a(v)$ appearing in the definition of the Cardy-Smirnov embedding. 
    {\bf Middle}: An embedding of a planar map into $\Delta$. 
    {\bf Right}: The percolation cycles of are given by the interfaces between blue and yellow vertices. In the scaling limit these converge to CLE$_6$.}
	\label{fig:cardy}
\end{figure}

The embedding considered in \cite{hs-cardy-embedding} is a novel discrete conformal embedding called the Cardy-Smirnov embedding, which is inspired by Smirnov's proof of Cardy's formula. Given a triangulation of a disk $M$ with three distinct boundary edges $a,b,c$ ordered counterclockwise,
we denote by $(a,b)$ the set of boundary vertices of $M$ between $a$ and $b$ in counterclockwise order (including one endpoint of $a$ and one endpoint of $b$).
Define $(b,c)$ and $(c,a)$ similarly. Let $\omega:\cV(M)\to\{ \op{blue},\op{yellow}  \}$ be a coloring of $M$ such that the boundary vertices are blue and the interior vertices are colored uniformly and independently at random in yellow and blue. For a vertex $v\in \cV(M)$, let
$E_a(v)$ be the event that there exists a simple path $P$ on $M$ such that 
\begin{compactitem}
	\item[(a)] $P$ has one endpoint in $(c,a)$ and one endpoint in $(a,b)$, while all other vertices of $P$ are inner blue vertices, and 
	\item[(b)] either $v\in P$ or $v$ is on the same side of $P$ as the edge $a$.
\end{compactitem}
See Figure~\ref{fig:cardy}. We define the events $E_b(v)$ and $E_c(v)$ similarly. 

Consider the equilateral triangle $\Delta :=\{(x,y,z):x+y+z=1,\,x,y,z>0 \}$ and let $\ol\Delta$ denote its closure. The Cardy-Smirnov embedding of $(M,a,b,c)$ is the function $\op{CdySmi}_M:\cV(M)\to\ol\Delta$ given by 
$$
\op{CdySmi}_M(v)= \frac{(\P_M [E_a(v)],\P_M[E_b(v)],\P_M[E_c(v)])}{\P_M [E_a(v)]+\P_M[E_b(v)]+\P_M[E_c(v)]},
$$ 
where $\P_M$ is the law of the percolation $\omega$ given the planar map $M$.

Smirnov showed that if we apply the Cardy-Smirnov embedding to the triangular lattice restricted to some simply connected Jordan domain $D$, then we approximate the Riemann mapping from $D$ to $\Delta$. Conformal invariance of percolation and Cardy's crossing formula are  immediate consequences of this result.

The embedded planar map defines an area measure $\cA_n$ and a (pseudo)metric $\frk D_n$ in $\ol\Delta$ as follows. If $M$ has perimeter $\ell$ we define $n=\ell^2$ and for any measurable set $U\subset\ol\Delta$ we define 
\eqb
\cA_n(U) = n^{-1}\cdot \# \{v\in \cV(M)\,:\,\op{CdySmi}_M(v)\in U  \}.
\label{eq:meas-emb}
\eqe
In other words, we give each embedded vertex mass $n^{-1}$. 
For $x\in \ol\Delta$ let $\frk v(x)\in \cV(M)$ be such that $\|x-\op{CdySmi}_M(\frk v(x))\|$ is minimized. Let $d_{\op{gr}}:\cV(M)\times \cV(M)\to\Z_+$ denote the graph distance and for $x,y\in\ol\Delta$ define
\eqb
\frk D_n(x,y) = n^{-1/4}\cdot d_{\op{gr}}( \frk v(x),\frk v(y) ).
\label{eq:metric-emb}
\eqe

Let $\cA$ denote the area measure of a quantum disk with perimeter 1 sampled from $\op{QD}(1)$ with $\gamma=\sqrt{8/3}$, embedded into $\Delta$ such that three uniform points from the LQG boundary measure $\cL$ are sent to the three vertices of $\Delta$. Let $\frk D$ denote the associated $\sqrt{8/3}$-LQG metric. The following theorem gives convergence of uniform triangulations of a disk under the Cardy-Smirnov embedding, where we use $\Rightarrow$ to denote convergence in law.
\begin{theorem}
	There exist deterministic constants $c_1,c_2>0$ such that $(c_1\cA_n,c_2\frk D_n)\Rightarrow (\cA,\frk D)$, where we equip the first coordinate with the weak topology and the second coordinate with the uniform topology.
	\label{thm:cardy}
\end{theorem}
We remark that one can obtain convergence of uniform triangulations with other topologies (e.g.\ whole-plane or sphere) by modifying the definition of $\op{CdySmi}_M$ appropriately. Furthermore, Theorem \ref{thm:cardy} also holds for  triangulations where we allow self-loops or do not allow multiple edges since these can be coupled to type II triangulations.

The key idea of the proof is to establish quenched convergence of site percolation on the uniform triangulation.  In order to describe what we prove, let $M$ be as in Theorem \ref{thm:cardy}, i.e., it is a type II Boltzman triangulation of a disk with fixed perimeter $\ell$. With Albenque \cite{aasw-type2}, building on \cite{legall-uniqueness,miermont-brownian-map,ps-triangulation,bet-mier-disk,aw-core,ab-simple}, we prove that $M$ converges in the scaling limit as a metric measure space (i.e., for the Gromov-Hausdorff-Prokhorov topology) to the Brownian disk with unit perimeter. It is known since \cite{lqg-tbm3} that this Brownian disk has a canonical conformal structure which makes it equivalent to a disk sampled from $\op{QD}(1)$. Let $a,b,c$ be three marked boundary edges chosen uniformly at random from $M$. For some $k\in\N$ consider independent percolations $\omega_1,\dots,\omega_k$ on $M$. Each percolation defines a loop ensemble (i.e., a countable collection of loops) as in the right panel of Figure \ref{fig:cardy}. We view $(M,\omega_1,\dots,\omega_k)$ as a metric measure space decorated by $k$ loop ensembles. We consider two loop ensembles $\omega$ and $\omega'$ to be close if, for any loop in $\omega$ (resp.\ $\omega'$) there is a loop in $\omega'$ (resp.\ $\omega$) which is close for the uniform topology modulo reparametrization of time. We call the resulting topology on tuples $(M,\omega_1,\dots,\omega_k)$ the \emph{Gromov-Hausdorff-Prokhorov-Loop}  topology. In the below theorem, $\cS$ is the metric measure space associated with a quantum disk for $\gamma=\sqrt{8/3}$ (equivalently, $\cS$ is the Brownian disk) and $\Gamma_1,\dots,\Gamma_k$ are $k$ independent CLE$_6$ on this disk.
\begin{theorem}\label{thm:joint}
	$(M,\omega_1,\dots,\omega_k)\Rightarrow (\cS,\Gamma_1,\dots,\Gamma_k)$ in the Gromov-Hausdorff-Prokhorov-Loop topology.	
\end{theorem}
Theorem \ref{thm:cardy} is an almost immediate consequence of the latter theorem: Crossing events in the discrete (resp.\ continuum) are encoded by the $\omega_j$'s (resp.\ the $\Gamma_j$'s). Using this and Theorem \ref{thm:joint}, we know that the conformal structure of $M$ can be studied by considering empirical crossing statistics for $\omega_1,\dots,\omega_k$ with $k$ large, since, by Theorem \ref{thm:joint}, these match the crossing statistics for $\Gamma_1,\dots,\Gamma_k$, which reflect the conformal structure of $\cS$.

Theorem \ref{thm:joint} can be viewed as a quenched convergence result for percolation on random planar maps.  Earlier annealed results were proved in \cite{gwynne-miller-perc,ghs-metric-peano} building on the mating-of-trees convergence of percolated triangulations (Theorem \ref{thm:perc}), which is also an annealed convergence statement for percolation on the planar map to CLE$_6$. In particular, \cite{ghs-metric-peano} written with Gwynne proves convergence of $(M,\omega_j)$ in the Gromov-Hausdorff-Prokhorov-Loop topology for each fixed $j\in\{1,\dots,k \}$. Building on this result, in order to prove Theorem \ref{thm:joint} it remains to prove independence of $\Gamma_1,\dots,\Gamma_k$ in any subsequential limit of $(M,\omega_1,\dots,\omega_k)$. The mating-of-trees approach alone is not sufficient for such convergence because the random walk only encodes the coupled information of the map and the statistical physics models, and the joint law of the mating-of-trees walk for 
$(M,\omega_1)$ and $(M,\omega_2)$ is hard to access. We achieve our goal  by studying the scaling limit of so-called \emph{dynamical percolation} on a random planar map. In dynamical percolation, each vertex of the random planar map has an independent Poisson clock (i.e., a random clock which rings at the set of times given by a Poisson point process) and the color of the vertex is resampled every time its clock rings. The continuum analogue of dynamical percolation is a dynamically changing CLE$_6$, and we call this process \emph{Liouville dynamical percolation} in our work with Garban and Sepulveda \cite{ghss-ldp}. The so-called pivotal points of the percolation play a key role for the process. The pivotal points are studied in detail in the works \cite{bhs-site-perc,hlls-cut-pts,hls-sle6} written with subsets of Bernardi, Lawler and Li. In particular, in the bijection from Theorem \ref{thm:perc} the pivotal points of $(M,\omega)$ are encoded in an explicit way by the walk $(\cL,\cR)$, which allows us to prove convergence of the pivotal measure. 
Independence of $\Gamma_1,\dots,\Gamma_k$ is deduced from a mixing result for Liouville dynamical percolation which is derived from the Euclidean counterpart due to Garban, Pete, and Schramm \cite{gps-fourier}.

%The proof proceeds by first proving convergence of the triangulated disk viewed as a metric measure space \cite{aasw-type2}. Then it is proved that this convergence is joint with a number of percolation observables on the planar map, including counting measure on pivotal points \cite{bhs-site-perc,ghs-metric-peano}. Building on these results, a mixing result \cite{ghss-ldp}, and a thorough understanding of the so-called pivotal points \cite{hlls-cut-pts,hls-sle6}, it is then argued that the scaling limit of dynamical percolation on the planar map is ergodic, which allows to conclude convergence of the embedded planar map since it implies quenched convergence of percolation on the planar map.

%The relevance of quenched scaling limits for convergence of planar maps under conformal embedding is also illustrated by the works \cite{}. In \cite{}, Gwynne, Miller, and Sheffield prove that 

Besides our work \cite{hs-cardy-embedding}, convergence under discrete conformal embeddings to LQG has been established for various random planar maps that can be viewed as coarse-grained versions of LQG surfaces.  One example is the so-called mated-CRT map. For each $\gamma\in (0,2)$, there exists a corresponding mated-CRT map that lies in the universality class of $\gamma$-LQG. They can be constructed in two ways: one from SLE-decorated decorated LQG surfaces and one from a pair of Brownian motions $(L,R)$. The equivalence of these two constructions follows from the mating-of-trees theorem~\cite{wedges}. Gwynne, Miller, and Sheffield~\cite{gms-tutte} proved that for $\gamma\in (0,2)$ the mated-CRT map converges to $\gamma$-LQG under the Tutte embedding. Another example is the Poisson-Voronoi tessellation of the Brownian disk, which converges to $\sqrt{8/3}$-LQG under the Tutte embedding \cite{gms-poisson-voronoi}. Both results are obtained by proving the quenched scaling limit for the simple random walk on the random planar map, which is based on the general result established in~\cite{gms-random-walk}. In these results only the area measure was shown to converge to the LQG area measure, while the metric was not treated.

A natural future research direction is to extend convergence results for planar maps under conformal embedding to other settings, for example other embeddings methods and other planar map models. Holden, partly in collaboration with Yu, are investigating these questions using mating-of-trees bijections, strong coupling between random walk and Brownian motion, and the method in \cite{gms-tutte}. Another question of major importance is metric convergence to LQG of non-uniform planar maps, which one can also attempt to investigate via the same strong couplings, strengthening the method from \cite{ghs-map-dist} mentioned above. Finally, it would be of significant interest to obtain scaling limit results for random planar maps without using the mating-of-trees framework. 

%established the conjecture in the special case of the Tutte embedding and the so-called mated-CRT map for $\gamma\in(0,2)$. The mated-CRT map is a particular random planar map obtained by discretizing a $\gamma$-LQG surface in a natural way. The Tutte embedding is one of the most commonly considered discrete conformal embeddings. It is defined for planar maps with disk topology by first embedding the boundary according to random walk hitting probabilities and then requiring that the embedding function is discrete harmonic at interior vertices; see Figure \ref{fig:tutte}. The authors first observe that if they have some other embedding $\wt\phi$ of the planar map such that random walk under this embedding converges to 2D Brownian motion (modulo reparametrization of time) then this embedding has to be close to the Tutte embedding, at least if $\wt\phi$ embeds the boundary vertices of the map in the correct cyclic order. They establish a general invariance principle for random walk in so-called ergodic scale-free environments satisfying a certain moment condition and connectivity property \cite{gms-random-walk} and argue that this applies in their setting. Variants of the strategy described here has also been applied to a few other settings \cite{gms-poisson-voronoi,afs-metric-ball,bgs-smith}.

\section{Liouville quantum gravity and conformal field theory.}
\label{sec:cft}

\newcommand{\rd}{\mathrm{d}}
\newcommand{\QA}{\mathrm{QA}}
\newcommand{\BA}{\mathrm{BA}}

\newcommand{\LF}{\mathrm{LF}}
\newcommand{\QS}{\mathrm{QS}}
\newcommand{\QD}{\mathrm{QD}}
\newcommand{\xin}[1]{{\color{blue}{#1}}}

A $d$-dimensional Euclidean \emph{quantum field theory} (QFT) is a way of assigning a measure on fields (i.e.\ functions) on all $d$-dimensional Riemann manifolds in a diffeomorphism invariant manner. The most basic example of a QFT is the Gaussian free field (GFF). When $d=2$, the GFF has the additional property that the field measure only changes by a constant multiple under conformal deformations. Namely, we have ${\mathrm{GFF}}[g](d\phi) =C_{\mathrm{Weyl}}[g,\omega]  {\mathrm{GFF}}[e^{\omega}g] (d\phi)$ where 
$g$ is a Riemannian metric on a manifold $M$, ${\mathrm{GFF}}[g]$ denotes GFF measure on $(M,g)$, 
$\omega$ is a real smooth function on $M$, and $C_{\mathrm{Weyl}}[g,\omega] $ is a constant. The constant multiple $C_{\mathrm{Weyl}}[g,\omega] $  has an explicit expression derived by Polyakov~\cite{polyakov-qg1} and is called the \emph{Weyl anomaly}. A 2D Euclidean QFT is called a \emph{conformal field theory} (CFT) if there exists a constant $\cc$ such that for each $(g,\omega)$ 
\begin{equation}\label{eq:wely}
 {\mathrm{Field}}[g](d\phi) = %C^{c}_{\mathrm{Weyl}}[g,\omega] 
 \big(C_{\mathrm{Weyl}}[g,\omega]\big)^{\cc}  {\mathrm{Field}}[e^{\omega}g] (d\phi),   
\end{equation}
where ${\mathrm{Field}}[g]$ denotes the law of the field. The exponent $\cc$ is called the \emph{central charge} of the CFT. 

In~\cite{polyakov-qg1}, Polyakov argued that 2D quantum gravity coupled with conformal matter can be described by the coupling of three CFTs: the matter CFT describing the conformal matter, the Liouville CFT describing the random geometry, and the ghost CFT describing the law of the conformal moduli. In Section~\ref{subsec:4.1} we review the rigorous construction of Liouville CFT  from~\cite{dkrv-lqg-sphere}. 
In Section~\ref{subsec:LCFT-MOT} we review the connection between Liouville CFT and the mating-of-trees framework, initiated in our joint work~\cite{ahs-sle,ahs-loop} with Ang. Here mating-of-trees refers to the coupling theory of SLE and LQG reviewed in Section 3. 
As we will see in Section~\ref{subsec:LCFT-MOT},  the rich interaction between Liouville CFT and mating of trees yields exact results on both sides which appear inaccessible from one side alone.
In Section~\ref{subsec:4.3}, we describe the matter-Liouville-ghost description and recent progress towards proving it.  Ideas reviewed in  Sections~\ref{subsec:LCFT-MOT} and~\ref{subsec:4.3} also have applications on the conformal matter side, namely SLE curves and scaling limit of 2D lattice models. In Section~\ref{subsec:4.4}, we review such applications using the example of percolation, which reveal intriguing CFT structure of this model.

\subsection{Background on Liouville conformal field theory.}\label{subsec:4.1} We first recall Liouville conformal field theory (LCFT) on the sphere, which is the QFT whose correlation functions are defined by the following path integral
\begin{align}\label{path_int_def}
\left \langle  \prod_{i=1}^N e^{\alpha_i \phi (z_i)} \right \rangle = 
\int_{\phi: \mathbb{H} \mapsto \mathbb{R}} D \phi    \prod_{i=1}^N e^{\alpha_i \phi(z_i)} e^{-S_{\gamma,\mu}[\phi;g]},
\end{align}
where $S_{\gamma,\mu}$ is the Liouville action given by
\begin{align}\label{liouville_action}
S_{\gamma,\mu}[\phi;g] = \frac{1}{4 \pi} \int_{\hat{\mathbb{C}}} \left( \vert \partial^{g} \phi \vert^2 + Q R_{g} \phi  \right  )
d v_{g} +  \mu  \int_{\hat{\mathbb{C}}}    e^{\gamma \phi} d v_{g}.   \end{align}
Here $\gamma\in (0,2)$, $Q=\frac{\gamma}{2}+\frac2{\gamma}$, and $g$ is a background metric on the sphere $\wh{\mathbb C}=\mathbb C\cup 
\{\infty\}$ with $R_{g}$ and $v_{g} $ being the curvature and volume, respectively. The constant $\mu>0$ is called the \emph{cosmological constant}. The rigorous construction of~\eqref{path_int_def} was achieved by David, Kupiainen, Rhodes, and Vargas~\cite{dkrv-lqg-sphere}. We  present this construction following the treatment of~\cite{ahs-sle}, which paves the way for connecting LCFT to mating of trees. 
The idea is that when $\mu=0$, the formal measure $D \phi \prod_{i=1}^N e^{\alpha_i \phi(z_i)}  e^{-S_{\gamma,0}[\phi;g]}$ can be realized via the GFF. For concreteness, we let $g$ be the metric such that $v_g=  |z|_+^4d^2z$  with $|z|_+:=\max\{|z|,1\}$. Then the measure $\LF_{\mathbb C}$ corresponding to $D \phi \, e^{-S_{\gamma,0}[\phi;g]}$ can be expressed by the \emph{whole-plane} GFF, which is  the variant of GFF with covariance kernel 
\begin{align}\label{eq:covariance}
G_\C(z,w) = -\log|z-w| + \log|z|_+ + \log|w|_+  \quad \textrm{for }z,w\in \C.\nonumber
\end{align}
\begin{definition}\label{def:LFc}
Let $P_{\mathbb C}$ be the law of whole-plane GFF above.	
Sample $(h, \mathbf c)$ from the measure $P_{\mathbb C} \times [e^{-2Qc}dc]$. Then the field  
 $\phi =  h -2Q \log |\cdot|_+ +\mathbf c$ is called the \emph{Liouville field on $\mathbb C$}. We denote its law by $\LF_{\mathbb C}$.
\end{definition}
The field measure $\LF_{ \C}^{(\alpha_i,z_i)_i}$ corresponding to $D \phi \prod_{i=1}^N e^{\alpha_i \phi(z_i)} e^{-S_{\gamma,0}[\phi;g]}$ can be realized as follows. 
\begin{definition}\label{def-RV-sph}
	Let $\alpha_i\in\mathbb{R}$  and $z_i\in \mathbb{C}$ for $i = 1, ..., N$ where $N\ge 1$ and all the $z_i$'s are distinct.
        Sample $(h, \mathbf c)$   from $ C _{\mathbb C}^{(\alpha_i,z_i)_i}  P_{\mathbb C} \times [e^{(\sum_i \alpha_i  - 2Q)c}dc]$ where  $C _{\mathbb C}^{(\alpha_i,z_i)_i}=\prod_{i=1}^N|z_i|_+^{-\alpha_i(2Q -\alpha_i)} e^{\sum_{i < j} \alpha_i \alpha_j G_{\mathbb C}(z_i, z_j)}$.
Then the field
\begin{equation}\label{eq:LF}
    \phi = h -2Q \log |\cdot|_+  + \sum_{i=1}^N \alpha_i G_{\mathbb C}(\cdot, z_i) + \mathbf c
\end{equation}
is called the  \emph{Liouville field on $ \mathbb C$ with insertions $(\alpha_i,z_i)_{1\le i\le N}$}. We denote its law by $\LF_{ \C}^{(\alpha_i,z_i)_i}$.
\end{definition}
Here the measures $\LF_{\mathbb C}$ and  $\LF_{\mathbb C}^{(\alpha_i,z_i)_i}$ are infinite because the free field measure $D \phi\,  e^{- \frac{1}{4 \pi} \int_{\hat{\mathbb{C}}}\vert \partial^{\hat g} \phi \vert^2  
d v_{g} }$ contains the zero mode $\mathbf c$. However, we still use terminologies in probability such as \emph{sample} and \emph{law} when dealing with infinite measures. Since samples from $\LF_{\mathbb C}^{(\alpha_i,z_i)_i}$ are variants of the GFF, it makes sense to define the Gaussian multiplicative chaos (GMC) measure $ e^{\gamma \phi} d^2z$. The term $-2Q \log |z|_+$ in~\eqref{eq:LF} ensures that $ \int_\C e^{\gamma \phi} d^2z$ is the correct  realization of ``$ \int_\C e^{\gamma \phi} dv_g$'' in the Liouville  action~\eqref{liouville_action}.
In conclusion, correlation functions for LCFT  are given  by
 \begin{equation}\label{eq:n-point}
\left \langle   \prod_{i=1}^N e^{\alpha_i  \phi (z_i)}   \right \rangle=  \int e^{-\mu \cA_\phi(\C)}  \LF_{ \C}^{(\alpha_i,z_i)_i}(d\phi ) 
\quad \textrm{where }  \cA_\phi\textrm{ is the GMC measure }  e^{\gamma \phi} d^2z\textrm{ on }\C.
\end{equation}

It was shown in~\cite{dkrv-lqg-sphere} that $\left \langle   \prod_{i=1}^N e^{\alpha_i  \phi (z_i)}   \right \rangle$ is finite and positive if and only if the following \emph{Seiberg bounds}  hold: 
\[\alpha_i<Q \textrm{ for } 1\le i\le N \quad \textrm{and}\quad \sum_{i=1}^{N}\alpha_i >2Q.\]
Furthermore, if the metric $g$ is changed to another conformally equivalent metric $e^{\omega}g$, the correlation functions change according to the Weyl anomaly~\eqref{eq:wely} with central charge $\cc=\ccL =1+6Q^2$. This means that LCFT is indeed a CFT with this central charge as argued by Polyakov~\cite{polyakov-qg1}. Finally, suppose $\psi:\mathbb{C}\to\mathbb{C}$ is a conformal map. Then 
\begin{equation}\label{eq:conf-map}
\left \langle   \prod_{i=1}^N e^{\alpha_i  \phi (\psi(z_i))}   \right \rangle=  \prod_{i=1}^N|\psi'(z_i)|^{-2\Delta_{\alpha_i}}   \left \langle   \prod_{i=1}^N e^{\alpha_i  \phi (z_i)}   \right \rangle,
\end{equation}
where $\Delta_\alpha=\frac{\alpha}{2}(Q-\frac{\alpha}{2})$ is called the  \emph{conformal weight}  of $e^{\alpha \phi(z)}$ in the CFT terminology. Consequently, the three-point function $\left \langle   \prod_{i=1}^3 e^{\alpha_i  \phi (\psi(z_i))}   \right \rangle$ must have the form
\begin{equation}\label{eq:DOZZ}
 C_{\gamma,\mu}(\alpha_1,\alpha_2,\alpha_3 )\times 
|z_1-z_2|^{2\Delta_{\alpha_3}-2\Delta_{\alpha_1}-2\Delta_{\alpha_2}}|z_1-z_3|^{2\Delta_{\alpha_2}-2\Delta_{\alpha_1}-2\Delta_{\alpha_3}}|z_2-z_3|^{2\Delta_{\alpha_1}-2\Delta_{\alpha_2}-2\Delta_{\alpha_3}},   
\end{equation}
where $C_{\gamma,\mu}(\alpha_1,\alpha_2,\alpha_3) $ is called the \emph{structure constant} of LCFT. 
In the groundbreaking work~\cite{krv-dozz}, Kupiainen, Rhodes, and Vargas derived an exact formula for $C_{\gamma,\mu}(\alpha_1,\alpha_2,\alpha_3) $,  which is the famous DOZZ formula predicted in the physics literature~\cite{do-dozz,zz-dozz}. In the axiomatic approach to CFT~\cite{bpz-conformal-symmetry}, the general correlation functions $\left \langle   \prod_{i=1}^N e^{\alpha_i  \phi (z_i)}   \right \rangle$ of a 2D CFT can be expressed in terms of the structure constant using a procedure called the \emph{conformal bootstrap}. This is rigorously established for LCFT by Guillarmou, Kupiainen, Rhodes, and Vargas~\cite{gkrv-sphere,gkrv-segal}. 

To connect LCFT to the mating-of-trees framework, it is crucial for us to consider LCFT on the disk. The rigorous construction is similar to the sphere case, as carried out in~\cite{hrv-disk}. Here we provide the counterpart of 
$\LF_{ \C}$ and $\LF_{ \C}^{(\alpha_i,z_i)_i}$.  The basic ingredient is the free-boundary GFF on the upper half plane $\bbH$, which is
the variant of GFF with covariance kernel 
\begin{align}\label{eq:covariance}
G_\bbH(z,w) = -\log |z-w| - \log|z-\ol w| + 2 \log|z|_+ + 2\log |w|_+  \quad \textrm{for } z,w\in \bbH\nonumber.
\end{align}
\begin{definition}\label{def:LFH}
Let $P_{\mathbb H}$ be the law of free-boundary GFF on $\bbH$.	Sample $(h, \mathbf c)$ from the measure $P_\mathbb{H}\times [e^{-Qc}dc]$. Then field  
 $\phi =  h -2Q \log |\cdot|_+ +\mathbf c$ is called the \emph{Liouville field on $\bbH$}. We denote its law by $\LF_{\bbH}$.
\end{definition}
\begin{definition}\label{def:LFH-ins}
Let $\beta_i\in\mathbb{R}$  and $s_i\in \partial\mathbb{H}$ for $i = 1, ..., M$ where $M\ge 1$ and all the $s_i$'s are distinct.
    Sample $(h, \mathbf c)$   from  $C_{\mathbb{H}}^{(\beta_i, s_i)_i}P_\mathbb{H}\times [e^{(\frac{1}{2}\sum_{i=1}^M\beta_i - Q)c}dc]$ with $C_{\mathbb{H}}^{(\beta_i, s_i)_i} =		\prod_{i=1}^M  |s_i|_+^{-\beta_i(Q-\frac{\beta_i}{2})} \exp(\frac{1}{4}\sum_{j=i+1}^{M}\beta_i\beta_j G_\mathbb{H}(s_i, s_j))$.   
Then  \[\phi = h - 2Q\log|\cdot|_++\frac{1}{2}\sum_{i=1}^M\beta_i G_\mathbb{H}(s_i, \cdot)+\mathbf{c}\] is called the  \emph{Liouville field on $ \mathbb H$ with insertions $(\beta_i,s_i)_{1\le i\le N}$}. We denote its law by $\LF_{ \bbH}^{(\beta_i,s_i)_i}$.
\end{definition}
Here $\LF_{ \bbH}^{(\beta_i,s_i)_i}$ is making sense of the field measure  $D \phi \prod_{i=1}^M e^{\frac{\beta_i}{2} \phi(s_i)} e^{-S_{\gamma,0}[\phi;g]}$ where $S_{\gamma,0}[\phi;g]$ is the Liouville action for $\phi:\bbH\mapsto \R$ in the absence of bulk and boundary cosmological constants. The correction function is defined by taking the joint Laplace transform of the total $\gamma$-LQG area of $\bbH$ and the $\gamma$-LQG lengths of line segments between boundary marked points. For example, the boundary three-point correlation function with $s_1<s_2<s_3$ is given by
 \begin{equation}\label{eq:3-point}
\left \langle   \prod_{i=1}^3 e^{\frac{\beta_j}2 \phi (s_i)}   \right \rangle=  \int e^{-\mu \cA_\phi(\bbH)-\sum_{i=1}^3 \mu_i \mathcal L_\phi(I_i)} \LF_{ \bbH}^{(\beta_i,s_i)_i}(d\phi),
 \end{equation}
where $\mu$ is the bulk cosmological constant, 
$\cA_\phi= e^{\gamma \phi} d^2z$ is the GMC measure on $\bbH$, 
$\mu_i$ ($1\le i\le 3)$ are the boundary cosmological constants, 
$\mathcal L_\phi=e^{\frac{\gamma}{2} \phi} dz$ is the GMC measure on $\partial\bbH$, and 
$I_2=[s_1,s_2]$, $I_3=[s_2,s_3]$, and $I_1=(-\infty,s_1]\cup [s_3,\infty)$.
Correlation functions with both bulk and boundary insertions can also be defined similarly by introducing the Liouville field 
$\LF_{ \bbH}^{(\alpha_i,z_i)_i, (\beta_j,s_j)_j}$. In particular,
for $z \in \mathbb{H}$ and $s \in \mathbb{R}$ we have
\begin{equation}\label{eq:1-1-point}
\left \langle  e^{\alpha \phi (z)}  e^{\frac{\beta}2 \phi (s)}   \right \rangle= \int e^{-\mu \cA_\phi(\bbH)- \mu_B \cL_{\phi}(\mathbb{R})}  \LF_{ \bbH}^{(\alpha,z), (\beta,s)}(d\phi),
 \end{equation}
 where $\mu$  and $\mu_B$ are the bulk and boundary cosmological constants, respectively. 
 We will not spell out the detailed definition of $\LF_{ \bbH}^{(\alpha,z), (\beta,s)}(d\phi)$ but refer to~\cite{ars-fzz}. 
 Again by conformal covariance,  $\left \langle \prod_{j=1}^3 e^{\frac{\beta_j}2 \phi (s_j)}    \right \rangle$ has the form
\begin{equation}\label{c4}
 \frac{H_\gamma(\beta_1,\beta_2,\beta_3;\mu,\mu_1,\mu_2,\mu_3)   }{\vert s_1 - s_2 \vert^{\Delta_{\beta_1} + \Delta_{\beta_2} - \Delta_{\beta_3}} \vert s_1 - s_3 \vert^{\Delta_{\beta_1} + \Delta_{\beta_3} - \Delta_{\beta_2}} \vert s_2 - s_3 \vert^{\Delta_{\beta_2} + \Delta_{\beta_3} - \Delta_{\beta_1}} },
\end{equation}  
while $\left \langle  e^{\alpha \phi (z)}  e^{\frac{\beta}2 \phi (s)}   \right \rangle$ has the form
\begin{equation}\label{c5}
\frac{G_\gamma(\alpha,\beta;\mu,\mu_B)}{\vert z - \overline{z} \vert^{2 \Delta_{\alpha} - \Delta_{\beta}} \vert z - s \vert^{2 \Delta_{\beta} }}.
\end{equation}  
In order to solve all the correlation functions of LCFT on the disk, one also needs to solve  $H_\gamma$ and $G_\gamma$ in~\eqref{c4} and~\eqref{c5} in addition to the three-point function on the sphere. They are called \emph{boundary structure constants}. The derivation of 
$H_\gamma$ and $G_\gamma$ requires the synergy between LCFT and the mating-of-trees framework. We defer this discussion to the end of Section~\ref{subsec:LCFT-MOT}.
LCFT on more general Riemannian manifolds were rigorously constructed in~\cite{drv-torus,remy-annulus,grv-higher-genus,Wu-bdy}. We refer to the  surveys~\cite{gkr-review,rv-icm} for further background on LCFT.

\subsection{Relation between LCFT and the mating-of-trees framework.}\label{subsec:LCFT-MOT}
Recall from Section 3 that the quantum sphere is the canonical LQG surface with  spherical topology that describes the scaling limit of random planar maps on the sphere. There are various versions of the quantum sphere depending on the number of marked points and on whether the area is fixed.  Let $\QS_n$ be the law of the free-area quantum sphere with $n$ marked points. 
In the mating-of-trees framework \cite{wedges}, one first defines $\QS_2$ in the cylindrical coordinate explicitly, and then defines $\QS_n$ for $n\neq 2$  by adding or removing marked points from $\QS_2$.  The fixed-area version  can be obtained by conditioning on the total area. It is a priori unclear how the quantum sphere and LCFT are related, despite sharing the same origin in physics~\cite{polyakov-qg1}.
The following theorem gives the precise link.
\begin{theorem}\label{def-QS-2}
Fix $\gamma\in (0,2)$. Let $(z_1, z_2, z_3) = (0, 1, e^{i\pi/3})$ and sample $\phi$ from $\frac{\pi \gamma}{2(Q-\gamma)^2}
\LF_\C^{(\gamma, z_1),(\gamma, z_2),(\gamma, z_3)}$. Then the law of the  quantum surface  $(\C, \phi,z_1, z_2, z_3)/{\sim_\gamma}$ is $\QS_3$.
\end{theorem}
The specific choice of $(z_1, z_2, z_3)$ in Theorem~\ref{def-QS-2} is only for concreteness. In fact, the argument for~\eqref{eq:conf-map} 
shows that  as $(z_1, z_2, z_3)$ vary, the law of $(\C, \phi,z_1, z_2, z_3)/{\sim_\gamma}$ only changes by a multiplicative constant.  
For the quantum disk, the analog of Theorem~\ref{def-QS-2} holds. Let $\QD_{n,m}$ be the law of the free-area-free-length quantum disk with $n$ interior marked points and $m$ boundary marked points. Sample $\phi$ from $\LF_\bbH^{(\gamma, s_1),(\gamma, s_2),(\gamma, s_3)}$ for some distinct $(s_1,s_2,s_3)\in \partial \bbH$. Then the law of the  quantum surface  $(\C, \phi,z_1, z_2, z_3)/{\sim_\gamma}$ is an explicit multiple of $\QD_{0,3}$. If $\LF_\bbH^{(\gamma, s_1),(\gamma, s_2),(\gamma, s_3)}$ is replaced by $\LF_\bbH^{(\gamma, z),(\gamma, s)}$ for some $z\in \bbH$ and $s\in \partial\bbH$, then the statement holds with $\QD_{0,3}$ replaced by $\QD_{1,1}$. 
The fixed-area version of Theorem~\ref{def-QS-2} was proved by Aru, Huang and Sun~\cite{ahs-sphere}, which implies Theorem~\ref{def-QS-2} without identifying the constant $\frac{\pi \gamma}{2(Q-\gamma)^2}$. In our joint work~\cite{ahs-sle} with Ang, Theorem~\ref{def-QS-2} was proved in a more conceptual manner with the explicit constant computed. A key insight from~\cite{ahs-sle} is that if we choose a random embedding of a sample from $\QS_0$ using the Haar measure on the M\"obius group (i.e. conformal automorphism of $\C$), 
then the law of the field on $\C$ is a multiple of $\LF_\C$. This random embedding is called the uniform embedding in~\cite{ahs-sle}.
The fixed-length version of the disk result was first proved by Cercle~\cite{cercle-quantum-disk}, following the method in~\cite{ahs-sphere}. Our work~\cite{ahs-sle} with Ang also applies to  $\QD_{0,3}$  with explicit constant, and the case of  $\QD_{1,1}$ was done in~\cite{ars-fzz} using the method from~\cite{ahs-sle}. 

The key mechanism that facilitates the rich interaction between  LCFT and mating of trees is the following meta statement first illustrated in our joint work~\cite{ahs-sle,ahs-loop} with Ang: LQG surfaces defined by  Liouville fields with insertions naturally fit into the mating-of-trees framework and behaves well under conformal welding. Here we provide two concrete instances which are the easiest to state without introducing more technical terms. 
They are immediate variants of the conformal welding results from~\cite{ahs-loop}.
\begin{theorem}\label{thm:QD-weld}
Let $(z_1, z_2, z_3) = (0, 1, e^{i\pi/3})$. Sample two independent copies of  $\QD_{1,1}$ with the restriction that their total boundary lengths agree. Conformally weld the two quantum surfaces together along the boundary with the boundary marked points identified. Then the law of the  resulting quantum surface agrees with $C\LF_\C^{(\alpha_1, z_1),(\alpha_2, z_2),(\alpha_3, z_3)} $ in the sense of~Theorem~\ref{def-QS-2}, where 
$(\alpha_1,\alpha_2,\alpha_3)=(\gamma,\gamma, \frac{\gamma}{2})$, 
$z_3$ lies on the interface of the two disks, and 
$C$ is an explicitly computable constant. 
Furthermore, the law of the interface is a variant of the SLE loop measure~\cite{zhan-loop-measures} with marked points. 
If we conformally weld two samples from  $\QD_{0,3}$, the analogous statement holds with all of $(z_1, z_2, z_3)$ being on the interface and $(\alpha_1,\alpha_2,\alpha_3)=(\frac{\gamma}{2}, \frac{\gamma}{2},\frac{\gamma}{2})$.
\end{theorem}
Nearby a point with a $\frac{\gamma}{2}$-insertion, the quantum surface locally look likes the neighborhood of the marked point on a two-pointed quantum sphere with weight $4$. Therefore  Theorem~\ref{thm:QD-weld} is locally compatible with the spherical version of Theorem~\ref{thm1}: the conformal welding of two copies of the two-pointed quantum disk of weight 2 is a two-pointed quantum sphere  of weight $4$. Starting from~\cite{ahs-sle,ahs-loop}, conformal welding results involving quantum surfaces described by LCFT have been established in various scenarios. For example, it was shown in~\cite{ahs-sle} that the two-pointed quantum disk with a generic weight also admits a natural LCFT description. It was shown in~\cite{asy-triangle} that when conformally welding  a two-pointed quantum disk to a surface arising from a Liouville field with three boundary insertions along one boundary arc, the resulting surface is described by a Liouville field with the three boundary insertions being modified accordingly. See~\cite{ay-radial,sy-simple,ahsy,hl-cle-nesting} for more instances of such conformal welding statements.

Both the mating-of-trees Brownian motion and  the  Liouville correlation functions encode the laws of quantum lengths and areas of certain quantum surfaces. Therefore, once both are put in the same playground via theorems such as Theorems~\ref{def-QS-2} and~\ref{thm:QD-weld},  one can hope to get exact results on both sides which appear inaccessible on one side.  
We give one example in each direction. On the one hand, Ang, Remy, and Sun derived the variance of the Brownian motions in the mating-of-trees framework~\cite{wedges}. 
On the other hand, Ang, Remy, Sun, and Zhu proved that the boundary structure constants $H_\gamma(\beta_1,\beta_2,\beta_3;\mu,\mu_1,\mu_2,\mu_3)$ and $G_\gamma(\alpha,\beta; \mu,\mu_B)$ 
for LCFT as defined in \eqref{c4} and~\eqref{c5} are given by the exact formulae predicted in physics~\cite{PT-bdy,Hosomichi}.   
The history behind these two results perfectly demonstrates the cross fertilization of the two sides. First, Remy and Zhu~\cite{remy-annulus,rz-gmc-interval,rz-bdy} derived the exact formulae for $H_\gamma$ and $G_\gamma$ 
 when $\mu=0$ using the proof method of~\cite{krv-dozz} for the DOZZ formula. The $\mu=0$ result was then used in~\cite{ars-fzz} to get the mating-of-trees variance, after which the formula for  $G_\gamma$ with $\mu>0$ and $\beta=\gamma$ was derived in the same paper using the SLE/LQG coupling. Finally, the complete formulae for $H_\gamma$  and $G_\gamma$ were established in~\cite{arsz} by a combination of the DOZZ argument and the SLE/LQG coupling.

 \subsection{The matter-Liouville-ghost CFT description for 2D quantum gravity.}\label{subsec:4.3}
 We start by recalling Polyakov's Liouville-ghost description of 2D pure gravity. 
 In theoretical physics, pure 2D quantum gravity on a topological surface $\mathcal S$ corresponds to a ``uniform'' measure on all possible 2D geometries on $\mathcal S$.  
Polyakov's seminal approach~\cite{polyakov-qg1} is to assume that random surfaces sampled from pure 2D quantum gravity are conformally equivalent to 2D smooth Riemann manifolds of the same topology, so that when the surface 
is conformally embedded onto a fixed Riemannian manifold, the random geometry is captured by the conformal factor. Therefore, the ``uniform'' measure on the space of all geometries on $\mathcal S$ can be expressed by the law of  $(\tau,\phi)$ where $\tau$ is a conformal structure on $\mathcal S$ and $\phi$ is the conformal factor. A key insight
from~\cite{polyakov-qg1}  is the proposal that conditioning on the Riemannian manifold, the law of the conformal factor  is governed by LCFT.  Namely the law of $(\tau,\phi)$ should have the form 
$C(g_\tau)\LF_{\mathcal S,g_\tau} (d\phi) D \tau$ where $D\tau$ is a natural volume measure on the (conformal) moduli space of  $\mathcal S$, $g_\tau$ is a smooth metric on $\mathcal S$ such that the Riemannian manifold $(\mathcal S,g_\tau)$ has  conformal structure $\tau$, $\LF_{\mathcal S,g_\tau}$ is the field measure given by LCFT on $(\mathcal S,g_\tau)$, and $C(g_\tau)$ is a $g_\tau$-dependent multiplicative factor. Based on some formal differential geometric reasoning, Polyakov argued  that $C(g_\tau)$  must behave as the partition of a CFT with central charge $-26$, which is now known as the b-c ghost  system~\cite{dms-cft-book}. Writing the ghost CFT partition function as $\mathcal Z_{\textrm{ghost}}[g_\tau]$, the law of $(\tau,\phi)$ can  then be written as
\begin{equation}\label{eq:mlg-BA}
\LF_{\mathcal S,g_\tau} (d\phi)\times \mathcal Z_{\textrm{ghost}}[g_\tau] D \tau.
\end{equation} 
Since the measure~\eqref{eq:mlg-BA} should not depend on choice of $g_\tau$, the Weyl anomaly from $\LF_{\mathcal S,g_\tau} $ and $\mathcal Z_{\textrm{ghost}}[g_\tau] $ must cancel. Therefore the LCFT describing the pure gravity must have central charge $\ccL=26$. Since  the $\gamma$-parameter and the central charge of LCFT are related by $\ccL=1+6Q^2$, we have $\gamma=\sqrt{8/3}$.

Uniform triangulation and quadrangulations are natural discrete models for pure 2D gravity. 
Their metric-measure scaling limits are given by \emph{Brownian surfaces}, as first shown for the sphere case by Le Gall~\cite{legall-uniqueness} and Miermont~\cite{miermont-brownian-map}. The Brownian disk and Brownian surfaces with  general topology were  introduced by Bettinelli and Miermont~\cite{bet-mier-disk,bet-mier-compact2}. 
%Therefore, the metric-measure structure of 2D pure gravity is captured by Brownian surfaces. 
As reviewed in Section~3, the  $\sqrt{8/3}$-LQG  sphere and disks describe the Brownian sphere and disk under conformal embeddings, respectively. By Theorem~\ref{def-QS-2} and its disk counterpart, we see that the law of the conformal factor under the conformal embedding is indeed governed by LCFT with $\gamma=\sqrt{8/3}$. Since the moduli spaces for sphere and disk are trivial, the ghost CFT does not enter the picture. Therefore, we conclude that for sphere and disk, Polyakov's CFT description of pure quantum gravity is completely verified.

The ghost part of Polyakov's description  for Brownian surfaces only manifests itself for non-simply connected surfaces. 
In the probabilistic setting, this was first formulated precisely by David, Rhodes, and Vargas~\cite{drv-torus} for the torus case, then by Guillarmou, Rhodes, and Vargas~\cite{grv-higher-genus} for higher genus closed surfaces, and by Remy~\cite{remy-annulus} for the annulus. The annulus case was proved by Ang, Remy and Sun~\cite{ars-annulus},
as an application of the interaction between LCFT and mating of trees. A crucial input for~\cite{ars-annulus} is the conformal bootstrap for LCFT on the annulus established by Wu~\cite{Wu-annulus}.  For concreteness, 
let $\mathcal A$ be the topological annulus and $g_\tau$ be such that $(\mathcal A,g_\tau)$ is the finite cylinder 
of unit circumference and length $\tau$. Then $Z_{\textrm{ghost}}[g_\tau]=\eta(2i\tau)^2$
where $\eta(z) = e^{\frac{i \pi z}{12}} \prod_{n=1}^\infty (1-e^{2 n i \pi z})$ is the Dedekind eta function.
\begin{theorem}\label{thm:ars-annulus}
Sample $(\phi,\tau)$ from $\LF_{\mathcal A,g_\tau} (d\phi) \times \mathcal Z_{\op{ghost}}[g_\tau] d\tau$, where $\LF_{\mathcal A,g_\tau}$ is the field measure for LCFT on $(\mathcal A,g_\tau)$ with $\gamma=\sqrt{8/3}$. Then  the law of the metric-measure space on $\mathcal S$ induced by $e^{\gamma \phi}$ agrees with the Brownian annulus coming from the scaling limit of uniform triangulations of annular topology~\cite{LeGall-annulus}. 
\end{theorem}
Proving  the Liouville-ghost description for Brownian surfaces other than the annulus remains an open question.
However, recently Baverez and Jego~\cite{bj25} gave a way of understanding the ghost CFT from conformal welding.  In a joint project with Baverez, Jego, Sun, and Wu  plan to prove the Liouville-ghost description for general Brownian surfaces based on~\cite{bj25}.

The CFT description of 2D quantum gravity coupled  with a conformal matter is given by 
\begin{equation}\label{eq:mlg-QA}
 \mathcal Z_{\op{M}}(g_\tau)\times \LF_{\mathcal S,g_\tau} (d\phi)\times \mathcal Z_{\textrm{ghost}}[g_\tau] D \tau,
\end{equation} 
where $\mathcal Z_{\op{M}}(g_\tau)$ is the partition function for the CFT describing the conformal matter. The cancellation of Weyl anomaly yields that the central charge $\ccL=1+6Q^2$ of LCFT is related to the central charge  $\ccM$ of the matter CFT  by $\ccL+\ccM=26$.
The decomposition~\eqref{eq:mlg-QA} was first argued by Polyakov~\cite{polyakov-qg1} when the matter is a free boson. This is the so-called \emph{bosonic string theory}. Guillamou, Rhodes, and Vargas~\cite{grv-higher-genus} showed that  the bosonic string theory is convergent when 
the free boson has central charge $\ccM\le 1$. When $\ccM=1$,  we have $\ccL=25$ and $\gamma=2$ for the LCFT, whose construction requires the critical Gaussian multiplicative chaos.

Besides free bosons, one can also consider 2D quantum gravity coupled with a general conformal matter  with $\ccM\le 1$, as done by David~\cite{david-conformal-gauge} and Distler-Kawai~\cite{dk-qg}. For example, when the conformal matter is a Virasoro minimal model from~\cite{bpz-conformal-symmetry}, one gets the minimal string theory~\cite{Seiberg-ministr}. 
Conformal matter with $\ccM\le 1$ can also be described  by 2D statistical physics models (e.g. the Ising model) at criticality. 
Therefore, random planar map decorated by statistical physics models give natural discrete models for 2D quantum gravity coupled with such conformal matters. Since SLE also describes the scaling limit of 2D statistical physics model at criticality, one can view SLE as an alternative approach to  conformal matters.
The relation between $\ccM$ and the $\kappa$-parameter for SLE is $\ccM=1-6(\frac{2}{\sqrt{\kappa}}-\frac{\sqrt{\kappa}}{2})^2$. Therefore $\ccM+\ccL=26$ yields  $\kappa=\gamma^2$ or $\kappa=16/\gamma^2$, which is consistent with the parameter relation in the SLE/LQG coupling in the mating of trees.  
In this sense, for sphere and disk where the ghost CFT is absent, the statement that the SLE/LQG coupling  describes the scaling limit of random planar map decorated by statistical physics models can be viewed as a version of the matter-Liouville decomposition of 
2D quantum gravity coupled with conformal matters with $\ccM\le 1$.  

The SLE/LQG perspective greatly enlarges the class of 2D quantum gravity coupled with conformal matters that can be analyzed quantitatively, even when the CFT description of the conformal matter is not fully understood.
For example, Ang, Remy and Sun~\cite{ars-annulus} considered quantum gravity  on the annulus with the conformal matter given by the conformal loop ensemble (CLE) or the SLE$_{8/3}$ loop measure~\cite{Werner-SAW}. The same approach for proving Theorem~\ref{thm:ars-annulus} also yields~\eqref{eq:mlg-QA} in these cases.  As we will see in Section~\ref{subsec:4.4}, this perspective leads to interesting applications to scaling limits of lattices models.  It would be interesting  to develop the mating-of-trees framework for general surfaces and link it to~the matter-Liouville-ghost decomposition.

\subsection{Applications to 2D statistical physics.}\label{subsec:4.4}
One of the key motivations for the mathematical understanding of Liouville quantum gravity is the application to scaling limits of 2D lattice models, which are conformally invariant random fractals. The Knizhnik-Polyakov-Zamolodchikov (KPZ) relation~\cite{kpz-scaling} is an explicit quadratic relation between the Euclidean scaling dimensions/exponents of these models and their counterpart in quantum gravity. The dimensions/exponents at the quantum level are sometimes possible to compute, and the KPZ relation yields the more interesting Euclidean results.  As a remarkable application of this idea, Duplantier~\cite{Duplantier1998} gave a convincing argument for the Mandelbrot 4/3-conjecture for planar Brownian motion by viewing it as the conformal matter. This conjecture was rigorously proved using SLE~\cite{lsw-bm-exponents3}.
Several rigorous versions of the KPZ relation have been established; see e.g.~\cite{shef-kpz,grv-kpz,gp-kpz}. In particular, Gywnne, Miller, and Holden~\cite{ghm-kpz} established a KPZ relation in the mating-of-trees framework and use it to derive dimension results for SLE curves. 

Inspired by the KPZ relation, in the joint work~\cite{ahs-sle} with Ang,  we demonstrated  that the synergy between mating of trees and LCFT can be used to derive exact results for SLE. The philosophy is the following. In light of the material reviewed in Section~\ref{subsec:LCFT-MOT}, one can quantitatively analyze a large family of geometric observables in the independent coupling of SLE and LCFT using the mating-of-trees framework. On the other hand, LCFT is exactly solvable thanks to the DOZZ formula, its boundary analogs, and the conformal bootstrap.
Therefore, one can extract exact solvability of SLE from the SLE/LCFT coupling. For 2D lattice models whose connection to SLE is established, we can use this method to get exact results on their scaling limits. We now use  percolation as an example to demonstrate that the SLE/LCFT coupling can yield results that appear inaccessible without employing the exact solvability of LCFT. Some results are even new in physics. We keep the percolation background minimal and refer to the introductory sections of~\cite{NQSZ2023,SXZ-Annulus,acsw-im-dozz} for further discussion and references.

Consider the critical Bernoulli site percolation on the triangular lattice where each vertex is colored black and white independently with probability  $\frac12$. Then its scaling limit can be described by SLE$_6$ and CLE$_6$ by the work of Smirnov~\cite{smirnov-cardy} and Camia-Newman~\cite{camia-newman-sle6}. For $R>r>0$, let $\mathbb A(r,R)=\{  r<|z|<R \}$ be the annulus with inner radius $r$ and outer radius $R$.  Then the probability that there exists a black path connecting the inner and outer boundary of $\mathbb A(r,R)$ converges to a limiting probability  $p_{B}(r,R)$ as the lattice  mesh size tends to 0.  We can similarly consider the event that there exist both a black path and a white path crossing the annulus 
$\mathbb A(r,R)$, and define the limiting crossing probability $p_{BW}(r,R)$. 
It was shown by Lawler, Schramm, and Werner that $p_{B}(r,R)\approx (r/R)^{-\alpha_1}$ with $\alpha_1=\frac{5}{48}$, and by  Smirnov and Werner~\cite{smirnov-werner-percolation} that  $p_{BW}(r,R)\approx (r/R)^{-\alpha_2}$ with $\alpha_2=\frac14$. It is equally natural to consider $p_{BB}(r,R)$, which describes the limiting probability for the event that there are two disjoint black paths crossing $\mathbb A(r,R)$. Then $p_{BW}(r,R)\approx (r/R)^{-\beta_2}$ where $\beta_2$ is known as the \emph{backbone exponent}. 
The value of $\beta_2$ was not known until its recent derivation by Nolin, Qian, Zhuang, and Sun~\cite{NQSZ2023} based on the SLE/LCFT coupling. They showed that $\beta_2$ is the unique solution in $(\frac{1}{4},\frac{2}{3})$ to the equation 
\begin{equation}\label{eq:backbone-exponent}
    \frac{\sqrt{36 x +3}}{4} + \sin \Big(\frac{2 \pi \sqrt{12 x +1}}{3} \Big) =0.   
\end{equation}
% \begin{theorem} [\cite{NQSZ2023}]
%    The backbone exponent $\beta_2$ is the unique solution in $(\frac{1}{4},\frac{2}{3})$ to 
%    \begin{equation}\label{eq:backbone-exponent}
%     \frac{\sqrt{36 x +3}}{4} + \sin \Big(\frac{2 \pi \sqrt{12 x +1}}{3} \Big) =0.   
% \end{equation}
% \end{theorem}
The key insight from~\cite{NQSZ2023} is that $\beta_2$ can be encoded in the law of the conformal radius of a random domain bounded by a certain SLE curve, which is defined by the derivative of a conformal map between the random domain to the unit disk.
Conformal derivatives naturally appear conformal coordinate change for LCFT correlation functions (e.g.~\eqref{eq:conf-map}), and their laws can be solved using the method from our previous work with Ang~\cite{ahs-sle}. A crucial ingredient for both~\cite{ahs-sle} and~\cite{NQSZ2023} is the boundary structure constants for LCFT derived by Remy and Zhu~\cite{rz-bdy}.
Similar ideas were used by Ang, Liu, Sun, Yu,  and Zhuang~\cite{ASYZ-touching,LSPZ-Fuzzy} to derive new scaling exponents for other percolation type models.

We conclude with a list of exact results for the conformal invariant scaling limits of 2D lattice models using the SLE/LCFT coupling. Sun, 
Xu, and Zhuang~\cite{SXZ-Annulus} obtained  exact formulae for  $p_{B}(r,R)$ and  $p_{BW}(r,R)$ predicted by Cardy~\cite{Cardy2002,Cardy2006} based on (non-rigorous) CFT consideration.  
Liu, Wu, and Zhuang~\cite{LWZ-mixing} derived the mixing rate exponent of 2D random-cluster from CLE, settling a  question of Duminil-Copin and Manolescu~\cite{DCM22}.
Cai, Fu, Sun, and Xie~\cite{CFSX2025} derived an exact formula for the probability that a planar Brownian motion disconnects the annulus $\mathbb A(r,R)$. 
Ang, Cai, Sun, and Wu~\cite{acsw-loop} derived the exact formula for the  electrical thickness of CLE loops. 
The same authors derived in~\cite{acsw-im-dozz} various three-point correlation functions for CLE on the sphere, including the formula conjectured by Delfino and Viti~\cite{DV-connectivity} for the connectivity of percolation cluster based on  CFT consideration.  The same quantity for the spin cluster of the three-Potts model was derived by Cai, Liu, Wu, and Zhuang~\cite{CLWZ25}, solving a problem from~\cite{Delfino_2013}. We expect the list to grow in the future. 

\section{Outlook.}\label{sec:outlook} We have already mentioned several open questions and future directions in Sections 3 and 4. We conclude this proceeding with a few additions.

\subsection{Connection to minimal string theory and JT gravity.} As mentioned in Section~\ref{subsec:4.3}, the minimal string theory is
given by 2D quantum gravity coupled with Virasoro minimal models. A non-perturbative solution to such string theories was discovered in physics around 1990 via the matrix model approach, which led to 
Witten's conjecture relating the intersection theory of algebraic curves to certain integrable systems~\cite{Witten}.  
Mirzakhani~\cite{Mirzakhani} gave a conceptual proof of  Witten's conjecture using a recursion  she derived for the Weil-Petersson volume of the moduli space of hyperbolic surfaces. Her recursion allows a matrix model approach to Jackiw--Teitelboim (JT)  gravity~\cite{SSS-JT}, a 2D gravity model of significant interest.  The minimal string theory is expected to be a deformation of JT gravity~\cite{LQG-JT}.  Both Virasoro minimal models and SLE describe the scaling limit of certain 2D lattice models, such as the Ising model. Therefore, SLE coupled with LQG provides an appealing geometric framework for minimal string theory as an alternative to matrix models. In particular, one can use SLE loops to cut LQG surfaces, similar in spirit to cutting hyperbolic surfaces by geodesic cycles, which is crucial in Mirzakhani's theory. It would be interesting to provide an exact encoding of Virasoro minimal models via SLE, construct the minimal string theory via the SLE/LQG coupling,  and study its relation to Mirzakhani's theory and JT gravity.

\subsection{CFT and 2D percolation.}
It is believed that geometrically interesting observables for critical 2D percolation are governed by a CFT. This yields profound predictions on the exact solvability of this model, including Cardy's formula and several results mentioned in~Section~\ref{subsec:4.4}.
Despite the substantial progress in proving such predictions, the precise nature of the CFT remains elusive, even at the physics level. 
A conformal bootstrap approach to percolation (and more general 2D random cluster and loop models) is being developed in physics. So far a number of three-point functions have been determined  at the physics level of rigor in terms of the so-called imaginary DOZZ formula. See \cite{Rilbault-survey} for a review. Recently, a relation between the imaginary DOZZ formula and CLE was established by Ang, Cai, Sun, and Wu~\cite{acsw-im-dozz}. Baverez-Jego~\cite{BJ24} and Gordina-Qian-Wang~\cite{GQW} constructed  Virasoro representation for the SLE loop measure. Guillarmou-Kupiainen-Rhodes~\cite{GKR-imLiouville}, Chatterjee~\cite{Chatterjee-imLiou} and  
Usciati-Guillarmou-Rhodes-Santachiara~\cite{UGRS} proposed path integral constructions for variants of LCFT whose three-point functions can be expressed by the imaginary DOZZ formula. Camia and Feng~\cite{Camia-Feng} revealed a logarithmic structure in percolation predicted in the conformal bootstrap framework. All these developments together with those reviewed in Section~\ref{subsec:4.4} give us hope that a complete understanding of the CFT structure underlying 2D percolation can be achieved one day.

\subsection{Random surfaces from Yang-Mills theory.} Polyakov's motivation for studying random surfaces comes from quantum gauge theory~\cite{Polyakov-book}, including the Yang-Mills theory.  By the so-called gauge/string duality, correlation functions in quantum gauge theory can  be expressed as a summation over surfaces. This was made precise for the lattice Yang-Mills theory by Chatterjee~\cite{Chatterjee-YM} and Cao-Park-Sheffield~\cite{CPS-YM}. One strategy to construct and study the continuous Yang-Mills theory is via the continuum limit of such random surfaces. The fundamental questions of mass gap and quark confinement can be naturally formulated in this approach. 
Random surfaces in this context are very different from the LQG surfaces reviewed in the current proceeding. In particular, these random surfaces need to be embedded into the $d$-dimensional Euclidean space  where the most interesting dimension is $d=4$. Moreover, the summation over surfaces involves negative weights and is absolutely divergent in the most interesting regime. 
Recently, Gwynne and Ding~\cite{Gwynne-Ding-super} constructed the random metric for LQG with central charge $\ccL\in (1,25)$. In addition,  Ang and Gwynne~\cite{Ang-Gwynne-mis,Ang-Gwynne-super} established a Markov property in the SLE/LQG coupling when the two parameters do not match as in the mating-of-trees framework. 
These developments provide a path to construct and study LQG surfaces embedded in $4$-dimensional Euclidean   space. We cherish the hope that this will inspire progress on Yang-Mills theory.

\section*{Acknowledgments.} We thank our advisor Scott Sheffield for introducing us the topic of Liouville quantum gravity and for his guidance and encouragement. We thank our collaborators on this topic, Marie Albenque, Morris Ang, Juhan Aru, Nathanael Berestycki, Oliver Bernardi, Jacopo Borga, Gefei Cai, Julien Dubédat, Hugo Falconet,  Xuesong Fu, Christophe Garban,  Ewain Gwynne,  Promit Ghosal, Yichao Huang, Greg Lawler, Matthis Lehmk\"uhler, Xinyi Li, Yiting Li, Haoyu Liu, Jason Miller, Pierre Nolin, Joshua Pfeffer, Ellen Powell, Wei Qian, Guillaume Remy, Scott Sheffield,  Avelio Sepulveda, Yi Sun, Sam Watson,  Baojun Wu, Zhuoyan Xie, Shengjing Xu, Pu Yu, Tunan Zhu, and Zijie Zhuang. We thank Morris Ang for helpful comments on an earlier draft.

% SIAM recommends using BibTeX
% if using BibTeX
\bibliographystyle{siamplain}

\def\cprime{$'$}

\end{document}